\documentclass[reqno]{amsart}

\usepackage{latexsym,amssymb,amsmath, amsfonts}
\usepackage{amsthm}
\usepackage{amscd}
\usepackage[shortlabels]{enumitem}
\usepackage{a4wide}

\numberwithin{equation}{section}
\usepackage[usenames,dvipsnames,svgnames,table]{xcolor}
\usepackage{hyperref}
\hypersetup{colorlinks=true,
	linkcolor=blue,
	urlcolor=blue,
	citecolor=red} 
\newcommand{\R}{\mathbb{R}}

\renewcommand{\le}{\leqslant}
\renewcommand{\ge}{\geqslant}
\renewcommand{\leq}{\leqslant}
\renewcommand{\geq}{\geqslant}

\newcommand{\be}{\begin{equation}}
\newcommand{\en}{\end{equation}}
\newcommand{\ee}{\end{equation}}

\newcommand{\bt}{\begin{theorem}}
\newcommand{\et}{\end{theorem}}
\newcommand{\bp}{\begin{proof}}
\newcommand{\ep}{\end{proof}}
\newcommand{\bc}{\begin{cor}}
\newcommand{\ec}{\end{cor}}
\newcommand{\bl}{\begin{lemma}}
\newcommand{\el}{\end{lemma}}
\newcommand{\bprop}{\begin{prop}}
\newcommand{\eprop}{\end{prop}}

\newcommand{\N}{\mathbb{N}}

\newtheorem{theorem}{Theorem}[section]

\newtheorem{lemma}[theorem]{Lemma}

\newtheorem{proposition}[theorem]{Proposition}

\numberwithin{theorem}{section} \numberwithin{definition}{section}

\newcommand{\RNum}[1]{\uppercase\expandafter{\romannumeral #1\relax}}
\def\R{\mathbb{R}}
\def\N{\mathbb{N}}

\newcommand{\vertiii}[1]{{\left\vert\kern-0.25ex\left\vert\kern-0.25ex\left\vert #1 
		\right\vert\kern-0.25ex\right\vert\kern-0.25ex\right\vert}}

\theoremstyle{definition}

\author[E. Bustamante]{Eddye Bustamante}
\address{Universidad Nacional de Colombia Sede Medell\'in, carrera 65 No. 59A -110, post 50034, Medell\'in Colombia}
\email{eabusta0@unal.edu.co}
\author[J. Jim\'enez]{Jos\'e Jim\'enez Urrea}
\address{Universidad Nacional de Colombia Sede Medell\'in, carrera 65 No. 59A -110, post 50034, Medell\'in Colombia}
\email{jmjimene@unal.edu.co}

\author[A. Mu\~noz]{Alexander Mu\~noz}
\address{Universidade de S\~ao Paulo, Instituto de Matemáticas e Estatística (IME-USP), R. do Matão, 1010 - Butantã, São Paulo - SP, 05508-090, Bra\-zil}
\email{alexd@usp.br}
\thanks{}

\date{}
\title[decay vs regularity of solutions to the 2D mZK equation]{On decay and regularity of solutions to the modified 2D Zakharov-Kuznetsov equation}
\keywords{Weighted Sobolev spaces, polynomial decay, local well-posedness}

\begin{document}

\begin{abstract} 
This work is devoted to study the relation between regularity and  decay for solutions of the two-dimensional modified Zakharov-Kuznetsov equation in the anisotropic weighted Sobolev spaces $Z_{s,(r_1,r_2)}:=H^s(\R^2)\cap L^2((1+|x|^{2r_1}+|y|^{2r_2})dxdy)$. 
\end{abstract}

\maketitle

\section{Introduction}\label{intro}

In this article we consider the initial value problem (IVP) associated to the two dimensional modified Zakharov-Kuznetsov (mZK) equation,
\begin{align}
\left. \begin{array}{rl}
u_t+\partial_x^3 u+\partial_x \partial_y^2 u +u^2 \partial_x u &\hspace{-2mm}=0,\qquad\qquad (x,y)\in\mathbb R^2,\; t\in\mathbb R,\\
u(x,y,0)&\hspace{-2mm}=u_0(x,y).
\end{array} \right\}\label{mZK}
\end{align}

This equation is a bidimensional generalization of the Korteweg-de Vries (KdV) equation and serves as a mathematical model to describe the propagation of nonlinear ion-acoustic waves in magnetized plasma (see for instance \cite{ZK1974}).\\

In recent years, the well-posedness of the IVP (1.1) has been studied in the context
of the classical Sobolev spaces $H^s(\mathbb{R}^2)$. The Zakharov-Kuznetsov equation, even with more general non-linearities than those we are considering, has been studied, among others, in \cite{BL2003, Faminskii1995, GH2014, K2022, K2021, LP2009} and references therein.
In the two dimensional case, Ribaud and Vento in \cite{RV2012} obtained well-posedness results for $s>1/4$ relying solely on linear estimates for free solutions.  Kinoshita addressed the critical case $s=1/4$ in \cite{K2022}, and also established $C^3$-ill-posedness for $s<1/4$. \\

In this work, we focus on the IVP \eqref{mZK} within the context of weighted Sobolev spaces. In this regard, the pioneering work for nonlinear dispersive equations is due to T. Kato \cite{KATO} for the Korteweg-de Vries (KdV) equation 
\begin{equation}
    \label{Kdv} \partial_t u +\partial_x^3 u +u \partial_x u=0.
\end{equation} Subsequently, in \cite{FLP} and \cite{NAHASPONCE2}, local well-posedness was established in weighted Sobolev spaces (for non-integer indices) for the Cauchy problem associated to \eqref{Kdv}. 

In the case of the Zakharov-Kuznetsov equation, the Cauchy problem in weighted Sobolev spaces was initially studied in \cite{BJM} and \cite{FP2016} for dimension $d=2$, where the authors adapted techniques developed for the KdV equation in \cite{FLP}, \cite{KATO}, and \cite{NAHASPONCE}. Additionally, Linares, Pastor, and Drumond Silva in \cite{LPD2021} also considered the local well-posedness in weighted Sobolev spaces for the 2D symmetrized Zakharov-Kuznetsov equation, as part of their work on dispersive blow-up phenomena. More recently, Mendez and Ria\~no in \cite[Theorem 1.2]{arg} obtained a first local well-posedness result for arbitrary dimensions $d\ge2$ without pursuing minimal regularity.  

It is worth mentioning that it is possible that the ideas of Ribaud and Vento \cite{RV2012} may be extended to the context of well-posedness in weighted Sobolev spaces. Since our main interest in this work is not to pursue optimal local well-posedness in such spaces, we will consider $s>3/4$ coming from the classical Kenig-Ponce-Vega method \cite{KPV1993} involving maximal estimates. \\

Isaza, Linares and Ponce in \cite[Theorem 1.4]{ILP2013} and \cite[Corollary 1.6]{ILP2} showed that there exists an improvement in regularity when an excess of decay is exhibited by the solutions of the KdV equation.  More precisely, they
proved that if $u\in C(\mathbb{R} ;L^2(\mathbb{R}))$ is the global solution of the KdV equation obtained in \cite{KPV1996}, and there exists $\alpha>0$ such that in two different times $t_0$, $t_1 \in \mathbb{R}$ 

\begin{align*}
|x|^{\alpha} u(x,t_0), \ |x|^{\alpha} u(x,t_1)\in L^2(\mathbb{R}), 
\end{align*}
then 
\begin{align*}
u \in C(\mathbb{R}:H^{2\alpha}(\mathbb{R})).
\end{align*}

Later on, in \cite{MP2}, Mu\~noz and Pastor revisited \cite[Theorem 1.4]{ILP2013} and adapted the proof for the modified KdV equation obtaining a similar result for any $\alpha>0$ without relying on the ideas of propagation of regularity done in \cite[Corollary 1.6]{ILP2}. A related result was also obtained for some dissipative perturbations of the KdV in \cite[Section 4]{A1}. Recently, Linares and Ponce in \cite{ILW} carry on the study of this phenomena obtaining a similar result for the Intermediate Long Wave equation. 

In the context of propagation of decay Mendez and Ria\~no in \cite[Theorem 1.1]{arg} proved that localized decay on the initial data is propagated by the corresponding solution. The authors showed that 
the gain of extra localized regularity is determined by the magnitude of the decay $r$ with respect to their maximum gain of regularity $\lfloor2r\rfloor$. In particular in the integer case, a weight of size $r$ produces regularity of order up to $2r$. 

Inspired by the ideas of Isaza, Linares and Ponce in \cite{ILP2013}, the main objective of this work is to study the relation between regularity and  decay for solutions of the two-dimensional modified Zakharov-Kuznetsov equation in the weighted Sobolev spaces 
\begin{align*}
Z_{s,(r_1,r_2)}:=H^s(\R^2)\cap L^2((1+|x|^{2r_1}+|y|^{2r_2})dxdy). 
\end{align*}

Taking into account the linear change of variables performed by Grünrock and Herr in \cite{GH2014}
\begin{align*}
x':=\mu x+\lambda y,\quad y':=\mu x-\lambda y, \quad t':=t\quad \text{and} \quad v(x',y',t'):=u(x,y,t),
\end{align*}
where $\mu=4^{-1/3}$ and $\lambda=\sqrt 3 \mu$, we consider the symmetrized version of equation \eqref{mZK}, that is
\begin{align}
v_t+\partial_x^3 v+\partial_y^3 v+ 4^{-1/3} (v^2\partial_x v+ v^2\partial_y v)=0.\label{smZK}
\end{align}
Note, $u$ satisfies the equation in \eqref{mZK} if and only if $v$ satisfies equation \eqref{smZK}. In this way, we may consider the IVP
\begin{align}
\left. \begin{array}{rl}
v_t+\partial_x^3 v+\partial_y^3 v+ 4^{-1/3} (v^2\partial_x v+ v^2\partial_y v) &\hspace{-2mm}=0,\qquad\qquad (x,y)\in\mathbb R^2,\; t\in\mathbb R,\\
v(x,y,0)&\hspace{-2mm}=v_0(x,y),
\end{array} \right\}\label{IVP_smZK}
\end{align}
where $v_0(x',y'):=u_0(x,y)$, instead of IVP \eqref{mZK}.\\

Following the ideas of Fonseca and Pachón in \cite{FP2016}, and Bustamante, Jiménez and Mejía  in \cite{BJM}, we will first establish in Theorem \ref{LWPZ} (below) a well posedness result for the IVP \eqref{IVP_smZK}. From Theorem \ref{LWPZ} and \cite[Theorem 1.1]{arg} the relationship between the regularity index $s$ and the decay rate $(r_1,r_2)$ of the solutions of the equation \eqref{smZK} is expected to be $r\mapsto2\max\{r_1,r_2\}$. Nevertheless, the main result of this article extends the relation of Mendez and Ria\~no to fractional exponents $s,r_1,r_2$ (in the particular case $d=2$) in the sense that the relationship between the indices $s$ and $(r_1,r_2)$ obeys $r\mapsto 2\min\{r_1,r_2\}$. More precisely,

\begin{theorem}\label{1.1}
    Let $v\in C([0,T]; H^s(\R^2))$ be the solution of the IVP \eqref{IVP_smZK} provided by Theorem \ref{LWPHS} (below). Assume there are $t_0<t_1\in [0,T]$ such that $(1+|x|^{s/2+\alpha_1}+|y|^{s/2+\alpha_2})v(t_i)\in L^2(\R^2)$, $i=0,1$. Then, $v\in C([0,T];H^{s+2\alpha}(\R^2))$, $\alpha:=\min \{ \alpha_1, \alpha_2 \}$ in the following cases
    \begin{enumerate}
        \item $\frac{3}{4}<s<1$ and $0\leq \alpha \leq \frac{1-s}{2}$.
        \item $s>\frac{9}{4}$ and $\alpha \geq 0$.
    \end{enumerate}
\end{theorem}
Note we may also interpret Theorem \ref{1.1} as an ill-posedness result in the following sense: assume  $v_0$ belongs to  $Y_{s,\alpha_1,\alpha_2}:=\left(H^s(\R^2)\setminus H^{s^+}(\R^2)\right)\cap L^2((1+|x|^{2\alpha_1}+|y|^{2\alpha_2})dxdy))$ with $\min\{\alpha_1,\alpha_2\}>s/2$. Then the corresponding solution of the mZK equation does not persist in $Y_{s,\alpha_1,\alpha_2}$.\\

The gap $s\in[1,\frac{9}{4}]$ ought to be compared to that displayed in Theorems 1.2 and 1.3 of \cite{ILW} in the context of the intermediate long wave equation. The argument used to prove our main result also involves a weighted energy estimate that remains standard when $\alpha\le \frac{1}{2}$. When $\alpha>1/2$ the structure of the nonlinear part of the equation is degenerated on increments of one derivative per each step of length $1/2$ on the size of $\alpha$.

We also point out that when $s>1$ the proof of our main result relies strongly on the use of maximal estimates in the presence of weights (see Lemma \ref{remark1} below) to handle the propagation of degeneration coming from derivatives of the nonlinear expression $f\partial f$. It is known for the non-symmetrized group that the estimate  $$\|U(t)u_0\|_{L^p_xL^\infty_{y,t\in[0,T]}}\leq c\|u_0\|_{H^s}$$ holds whenever \begin{equation}
    \label{cond1}s>1-1/p \mbox{ if } p\geq 4 \mbox{ or when } s>3/4 \mbox{ if } p\in[2,4).
\end{equation}
Moreover, this is sharp up to the endpoint (see Theorem 1.2 in \cite{LR2021}).  To control the weights in $H^s(\mathbb{R}^2)$ with $s>1$ in our method, the condition \eqref{cond1} is translated into $s>9/4$. It is expected this is the best we can get with this method. \\

This article is organized as follows: in Section 2 we establish some linear estimates for the group associated to the linear part of the symmetrization of the mZK equation, we present the Leibniz rule for fractional derivatives,  and we recall some results concerning the Stein derivative and interpolation estimates. Section 3 is devoted to the study of the local theory for the Cauchy problem  \eqref{IVP_smZK} in the weighted Sobolev spaces $Z_{s,(r_1,r_2)}$. In Section 4 we present the main contributions of this paper, we establish the relationship between decay and regularity of solutions for the mZK equation and prove Theorem \ref{1.1}.

\section{Preliminaries}\label{presection}
\subsection{Notation}\hfill

We use $c$ to denote several constants that may vary from line to line. Dependence on parameters is indicated using parenthesis. We write $a\sim b$ when there exists a constant $c>0$ such that $a/c\le b\le ca$. For a real number $s$ we denote with $s^+$ (respec. $s^-$) the quantity $s+\varepsilon$ (respec. $s-\varepsilon)$ for an arbitrary small $\varepsilon>0$.

The set $L^p(\R)$ is the usual Lebesgue space and $L^2(wdx)$ represent the space $L^2$ with respect to the measure $w(x)dx$. The norm in $L^p(\R)$ will be denoted by $\|\cdot\|_p$. For a function $f$, the expressions $\widehat{f}$ and $f^\vee$ mean the Fourier and inverse Fourier transform, respectively.  For $s, r_1, r_2\in \R$ the set $H^s(\R^2)$ is the $L^2$-based Sobolev space of index $s$ and the set
$$Z_{s,(r_1,r_2)}:=H^s(\R^2)\cap L^2((1+|x|^{2r_1}+|y|^{2r_2})dxdy)$$
is the weighted Sobolev space with regularity $s$ and decay $(r_1,r_2)$. In case $r_1=r_2$ we simply write $Z_{s,(r_1,r_2)}=:Z_{s,r_1}$.
We denote with $\|\cdot\|_{L^p_xL^q_T}$ the norm in the mixed Lebesgue spaces defined by \begin{equation*}
    \|f\|_{L^p_xL^q_T}:= \left( \int_\mathbb{R} \left( \int_0^T |f(x,t)|^q dt\right)^{p/q} dx
\right)^{1/p}.
\end{equation*}
We use $(\xi,\eta)\in\R^2$ for the dual Fourier variable of $(x,y)\in\mathbb R^2$. We denote the fractional derivatives with
\begin{align*}
(D^sf)^\wedge(\xi,\eta)&=(D_{xy}^sf)^\wedge(\xi,\eta):=|(\xi,\eta)|^s \widehat{f}(\xi,\eta),\\
(D_x^{s} f)^\wedge (\xi,\eta)&:=|\xi|^{s} \widehat f (\xi,\eta),\\
(D_y^{s} f)^\wedge (\xi,\eta)&:=|\eta|^{s} \widehat f (\xi,\eta).
\end{align*}

For simplicity we adopt the notation $\langle z\rangle:=(1+z^2)^{1/2}$. Also, denote with $J^sf$ the potential
$$\widehat{J^s(f)}(\xi,\eta)=(1+x^2+y^2)^{s/2}\widehat{f}(\xi,\eta).$$

If $f_m\to f$ and we say an estimate is valid for $m$ large enough we mean that for some $M\in \N$ if $m\ge M$ then
$$\|f_m\|\le \|f_m-f\|+\|f\|\le 1+\|f\|\le c\|f\|.$$
\subsection{Linear estimates}\hfill

Let us consider the linear IVP associated to the IVP \eqref{IVP_smZK},
\begin{align}
\left. \begin{array}{rl}
v_t+\partial_x^3 v+\partial_y^3 v &\hspace{-2mm}=0,\qquad\qquad (x,y)\in\mathbb R^2,\; t\in\mathbb R,\\
v(x,y,0)&\hspace{-2mm}=v_0(x,y),
\end{array} \right\}\label{lsZK}
\end{align}

whose solution is given by
$$v(x,y,t)=[V(t)v_0](x,y),\quad (x,y)\in\mathbb R^2,\quad t\in\mathbb R,$$
where $\{V(t):t\in\mathbb R\}$ is the unitary group defined by
\begin{align}
[V(t)v_0](x,y):=\frac1{2\pi}\int_{\mathbb R^2} e^{i[t(\xi^3+\eta^3)+x\xi+y\eta]} \hat v_0(\xi,\eta) d\xi d\eta.\label{G}
\end{align}

Lemmas \ref{STE} to \ref{MTE}, as well as Proposition \ref{P4}, deal with some results concerning the group defined in \eqref{G}.
\begin{lemma}\label{STE} (Strichartz-type estimates). Let $\varepsilon\in[0,1/2]$, $\gamma=\frac{(1-\varepsilon)}6$ and $\delta=(2-3\varepsilon)/18$. We have that
\begin{equation}
    \|V(\cdot_t)f \|_{L^3_T L^\infty_{xy}}\leq c T^\gamma \|D_x^{-\varepsilon/2}f \|_{L^2_{xy}}\label{Stri1_x},
\end{equation}
\begin{equation}
    \|V(\cdot_t)f \|_{L^3_T L^\infty_{xy}}\leq cT^\gamma \|D_y^{-\varepsilon/2}f \|_{L^2_{xy}}\label{Stri1_y},
\end{equation}
\begin{equation}
    \|V(\cdot_t)f \|_{L^{9/4}_T L^\infty_{xy}}\leq cT^\delta \|D_x^{-\varepsilon/2}f \|_{L^2_{xy}}\label{Stri2_x}
\end{equation}
and
\begin{equation}
    \|V(\cdot_t)f \|_{L^{9/4}_T L^\infty_{xy}}\leq cT^\delta \|D_x^{-\varepsilon/2}f \|_{L^2_{xy}}\label{Stri2_y}.
\end{equation}
\end{lemma}

\begin{proof} It is clear that the procedure developed in proving Lemma 2.1 in \cite{BJM}, and Lemma 2.6 in \cite{LP2009} can be adapted to establish the estimates \eqref{Stri1_x} to \eqref{Stri2_y}.

\end{proof}

Although the following Lemma is not strictly demonstrated in \cite{BJM} or \cite{LP2009}, it can be established using the same methods as shown in those works.

\begin{lemma}\label{LTE} (Local-type estimates). There exists a constant $c$ such that

\begin{equation}
    \| \partial_x V(\cdot_t)v_0\|_{L_x^\infty L^2_{ty}}\leq c \|v_0 \|_{L^2_{xy}},\label{LTE_x}
\end{equation}
\begin{equation}
    \| \partial_y V(\cdot_t)v_0\|_{L_y^\infty L^2_{tx}}\leq c \|v_0 \|_{L^2_{xy}},\label{LTE_y}
\end{equation}
\begin{equation}
    \left\| \partial_x \int_0^t V(t-t') g(\cdot,t') dt' \right\|_{L^\infty_T L^2_{xy}}\leq c \|g\|_{L^1_x L^2_{yT}}\label{LTE_xd}
\end{equation}
and
\begin{equation}
    \left\| \partial_y \int_0^t V(t-t') g(\cdot,t') dt' \right\|_{L^\infty_T L^2_{xy}}\leq c \|g\|_{L^1_y L^2_{xT}}.\label{LTE_yd}
\end{equation}

\end{lemma}
\begin{proof} Notice that estimates \eqref{LTE_xd} and \eqref{LTE_yd} are simply the dual versions of \eqref{LTE_x} and \eqref{LTE_y}, respectively. For the proof of \eqref{LTE_x} and \eqref{LTE_y} we refer the reader to Lemma 2.2 in \cite{BJM}.
\end{proof}

\begin{lemma}\label{MTE} (Maximal-type estimates). Let $v_0\in H^s(\mathbb R^2)$, for some $s>3/4$. Then for all $T>0$
\begin{align}
\| V(\cdot_t)v_0 \|_{L_x^2 L_{yT}^\infty}\leq c(s)(1+T)^{1/2} \|v_0 \|_{H^s_{xy}},\label{MTE_x}\\
\| V(\cdot_t)v_0 \|_{L_y^2 L_{xT}^\infty}\leq c(s) (1+T)^{1/2} \|v_0 \|_{H^s_{xy}}. \label{MTE_y}
\end{align}
\end{lemma}
\begin{proof}
    See Lemma 2.3 in \cite{BJM}.
\end{proof}

\begin{proposition}\label{P4}
    For $6\le r\le\infty$ and $s=\frac{1}{2}-\frac{3}{r}$ we have
    \begin{equation}\label{grunrock}
        \|V(t)f\|_{L^4_{xy}L^r_t}\le c\|f\|_{\dot{H}_{xy}^s}.
    \end{equation}
\end{proposition}
\begin{proof}
    See Corollary 1 in \cite{Grua}.
\end{proof}

Using \eqref{grunrock} and Sobolev embeddings it can be seen that for $s>3/4$ it follows
\begin{equation}
    \label{grunrockcol} \|V(t)f\|_{L^4_xL^\infty_{yt}}\le c\|f\|_{H^s_{xy}} \mbox{ and }  \|V(t)f\|_{L^4_yL^\infty_{xt}}\le c\|f\|_{H^s_{xy}}.
\end{equation}

\subsection{Fractional derivatives}\hfill

The following inequality is crucial to deal with the $L^2-$norm of the fractional derivative of a product.

\begin{lemma}\label{LFR} (Leibniz fractional rule). Let us consider $0<\alpha<1$ and $1<p<\infty$. Thus
\begin{align*}
\| D^\alpha(fg)-f D^\alpha g-gD^\alpha f \|_{L^p(\mathbb R)}\leq c \|g \|_{L^\infty(\mathbb R)} \| D^\alpha f\|_{L^p(\mathbb R)},
\end{align*}
where $D^\alpha$ denotes $D^\alpha_x$ or $D^\alpha_y$.
\end{lemma}
\begin{proof}
    See Theorem A12 in \cite{KPV1993}.
\end{proof}

In what comes to the $L^p-$norm of the fractional derivative of a function, we have the following useful result. 
\begin{lemma}\label{Int_Der} Let $l,m$, $r\in(1,\infty]$, and $0<a<b$ with
$$\frac1l=\left(1-\frac ab)\right)\frac1r+\frac ab\frac1m.$$
Then
$$\|D^af\|_{L^p}\leq c \|f\|_{L^r}^{1-a/b}\|D^b f\|^{a/b}_{L^m},$$
where $D^\alpha$ denotes $D^\alpha_x$ or $D^\alpha_y$.
\end{lemma}
\begin{proof}
    See Lemma 5 in \cite{NAHASPONCE}.
\end{proof}
%
%

Denote with $A_p$ the Muckenhoupt class on $\R^n$. More precisely, given $1<p<\infty$, $A_p$ consists of all weights $\omega$ such that
\begin{equation}
	[\omega]_p=\sup_{Q}\left(\frac{1}{|Q|}\int_{Q}\omega(y)dy \right)\left(\frac{1}{|Q|}\int_{Q}\omega^{-\frac{1}{p-1}}(y)dy \right)^{p-1}<\infty,
\end{equation}
where the supremum is taken over all cubes $Q\subset\R^n$; additional details may be seen in \cite{Hunt}.

\begin{lemma}\label{hilberpp}
The Hilbert transform is bounded in $L^p(wdx)$, $1<p<\infty$, if and only if $w\in A_p$.\vspace{-2mm}
\end{lemma}
\begin{proof} See Theorem 9 in \cite{Hunt}.
\end{proof}
\subsection{Interpolation results}\hfill

\begin{lemma}\label{interplemma}
	Assume $a, b>0$, $1<p<\infty$ and $\theta\in(0,1)$. If $J^a f\in L^p(\R^n)$ and $\langle x \rangle^b f\in L^p(\R^n)$ then 
	\begin{equation}\label{interp2}
		\|\langle x\rangle^{\theta b}J^{(1-\theta) a}f\|_{L^p(\R^n)}\le c \|\langle x\rangle^b f\|_{L^p(\R^n)}^{\theta}\|J^a f\|_{L^p(\R^n)}^{1-\theta}.
	\end{equation}
	The same holds for $D$ instead of $J$. Moreover, for $p=2$ we have
	\begin{equation}\label{interp}
		\|J^{\theta a}\langle x\rangle^{(1-\theta)b}f\|_{L^2(\R^n)}\le c \|J^a f\|_{L^2(\R^n)}^\theta \|\langle x\rangle^b f\|_{L^2(\R^n)}^{1-\theta}.
	\end{equation}
\end{lemma}
\begin{proof}
	Inequality	\eqref{interp} follows from \eqref{interp2} in view of Plancherel's identity. The proof of \eqref{interp2} follows using Hadamard's three lines theorem. See Lemma 4 in \cite{NAHASPONCE}.
\end{proof}

Under slight modifications of the proof of Lemma \ref{interplemma} it can be seen that for all $j\in\{1,2,\dots,n\}$ it follows
\begin{equation}
    \label{interpol1r} \|\langle x_j\rangle^{\theta b}J^{(1-\theta) a}f\|_{L^2(\R^n)}\le c \|\langle x_j\rangle^b f\|_{L^2(\R^n)}^{\theta}\|J^a f\|_{L^2(\R^n)}^{1-\theta}
\end{equation}
and
\begin{equation}
    \label{interpol1} \|J^{\theta a}\langle x_j\rangle^{(1-\theta)b}f\|_{L^2(\R^n)}\le c \|J^a f\|_{L^2(\R^n)}^\theta \|\langle x_j\rangle^b f\|_{L^2(\R^n)}^{1-\theta}.
\end{equation}

Consider the space $L^2_{j,\alpha}(\R^n):=\{f\in L^2(\R^n) \mid (1+\xi_j^2)^{\alpha/2}\widehat{f}(\xi)\in L^2(\R^n) \}$ and define the directional Stein derivative \begin{equation}
    \label{steind} \mathcal{D}_j^{\alpha}(f)(x)=\left(\int_{\R} \frac{|f(x+t\Vec{e}_j)-f(x)|^2}{|t|^{1+2\alpha}} dt\right)^{1/2}.
\end{equation}

\begin{theorem}\label{rhcp}
    Let $\alpha\in(0,1)$. $f\in L^2_{j,\alpha}(\R^n)$ if and only if $f\in L^2(\R^n)$ and $\mathcal{D}_{j}^\alpha(f)\in L^2(\R^n)$. Moreover $\|J_{j}^{\alpha}f\|_{L^2} \sim \|f\|_{L^2}+\|\mathcal{D}_j^{\alpha}\|_{L^2}$.
\end{theorem}
\begin{proof}
    Note $$ \|\mathcal{D}_j^\alpha(f)\|_{L^2}^2=\int_\R \frac{\|f(\cdot+t\vec{e}_j)-f(\cdot)\|_{L^2_x}^2}{|t|^{1+2\alpha}}dt.$$
   Since \begin{equation*}
       \|f(\cdot+t\vec{e}_j)-f(\cdot)\|_{L^2_x}^2=\int_{\R^n}|e^{it\vec{e}_j\cdot \xi}-1||\widehat{f}(\xi)|^2d\xi
   \end{equation*}
   we have 
   \begin{equation}
       \begin{split}
           \|\mathcal{D}_j^\alpha(f)\|_{L^2}^2&=\int_{\R^n}|\widehat{f}(\xi)|^2\int_\R \frac{|e^{it\xi_j}-1|^2}{|t|^{1+2\alpha}}dtd\xi\\&=\int_{\R^n} |\widehat{f}(\xi)|^2 |\xi_j|^{2\alpha}\int_\R \frac{|e^{iu}-1|^2}{|u|^{1+2\alpha}}dud\xi=c(\alpha)\int_{\R^n}|\widehat{f}(\xi)|^2 |\xi_j|^{2\alpha}d\xi.
       \end{split}
   \end{equation}
We conclude that the fact that both $\mathcal{D}_j^\alpha(f)$ and $f$ are in $L^2(\R^n)$, is equivalent to $(1+\xi_j^2)^{\alpha/2}\widehat{f}\in L^2(\R^n)$. 
\end{proof}

\begin{lemma}\label{interp1}
    Let $a,b>0$. For $j\in\{1,\dots,n\}$ assume that $\widehat{J^a_{x_j}f}=(1+\xi_j^2)^{a/2}\widehat{f}$ and $\langle x_j\rangle^{b}f=(1+|x|^2)^{b/2}f$ are in $L^2(\R^n)$. Then, for all $\theta\in(0,1)$ we have
    \begin{equation}
        \|J^{a\theta}_{x_j} \langle x_j\rangle^{(1-\theta)b}f\|_{L^2}\le c\|J^a_{x_j} f\|_{L^2}^{\theta}\|\langle x_j\rangle^{b}f \|_{L^2}^{1-\theta}.
    \end{equation}
\end{lemma}
\begin{proof}
    It is enough to consider $a=1+\beta$, with $\beta\in(0,1)$. For any $g\in L^2(\R^n)$ define $$F_j(z):=e^{z^2-1}\int_{\R^n}J^{az}_{x_j}\left(\langle x_j\rangle^{(1-\theta)b}f(x)\right)\overline{g(x)}dx.$$
    Using Theorem \ref{rhcp} it can be seen, in a similar fashion as done in the proof of Lemma \ref{interplemma}, that $$|F_j(iw)|\le e^{-w^2-1}\|\langle x_j\rangle^{b}f\|$$ and $$|F_j(1+iw)|\le c(\alpha)e^{-w^2}(1+|wb|^2)\|J^{1+\beta}_{x_j}f\|_{L^2}.$$
    The lemma follows from Hadamard's three lines theorem and duality. 
\end{proof}

\section{Local theory}

We gather some useful results regarding the existence of solutions to the IVP \eqref{IVP_smZK} in Sobolev and Weighted spaces. We begin with the result concerning the local well-posedness in $H^s(\mathbb R^2)$.

\begin{theorem}\label{LWPHS}
Let $s>3/4$ and $v_0\in H^s(\mathbb R^2)$ a real valued function. Then there exist $T>0$ and a unique $v$, in a certain subspace $X_{T,s}$ of $C([0,T];H^s(\mathbb R^2))$, solution of the IVP \eqref{IVP_smZK}. (The definition of the subspace $X_{T,s}$ will be clear in the proof of the theorem), such that
\begin{equation*}
    \begin{split}
        &\| D^s_x v_x\|_{L^\infty_x L^2_{yT}} + \|D^s_y v_x\|_{L^\infty_x L^2_{yT}} +\|D^s_y v_x\|_{L^\infty_y L^2_{xT}} + \|D^s_x v_y\|_{L^\infty_y L^2_{xT}} < \infty,\\&
\|v\|_{L^3_T L^\infty_{xy}} + \|v_x\|_{L^{9/4}_T L^\infty_{xy}} \|v_y\|_{L^{9/4}_T L^\infty_{xy}}<\infty\ \mbox{ and  }\ \|v\|_{L^2_x L^\infty_{yT}} + \|v\|_{L^2_y L^\infty_{xT}}<\infty.
    \end{split}
\end{equation*}

Moreover, for any $T'\in(0,T)$ there exists a neighborhood $V$ of $v_0$ in $H^s(\mathbb R^2)$ such that the data-solution map $\tilde v_0 \mapsto \tilde v$ from $V$ into $X_{T'}$ is Lipschitz.
\end{theorem}

\begin{proof} The proof of this result employs the contraction principle and follows the ideas of Linares and Pastor in \cite{LP2009} for the modified (non-symmetrized) ZK equation; and of Bustamante, Jiménez, and Mejía in \cite{BJM} for the ZK equation after being symmetrized. Due to the similarity of the estimates involved, we will provide only a sketch of the proof.\\

We begin by considering the integral operator
\begin{align}
\Psi(v)(t):= V(t) v_0 - 4^{-1/3} \int_{-t}^t V(t-t')(v v_x + v v_y)(t') dt'\equiv V(t)v_0 + E(v)(t),\label{IntEq}
\end{align}
and defining the metric space
$$X_{T,s}:=\{v\in C([0,T]; H^s(\mathbb R^2)): \eta_s(v) <\infty\},$$

where
\begin{align}
\notag \eta_s(v):= & \|v\|_{L^\infty_T H^s_{xy}} + \| D^s_x v_x\|_{L^\infty_x L^2_{yT}} + \|D^s_y v_x\|_{L^\infty_x L^2_{yT}} + \|v_x\|_{L^{9/4}_T L^\infty_{xy}} + \|v\|_{L^2_x L^\infty_{yT}} + \|D^s_x v_y\|_{L^\infty_y L^2_{xT}}\\
& + \|D^s_y v_x\|_{L^\infty_y L^2_{xT}} + \|v_y\|_{L^{9/4}_T L^\infty_{xy}} + \|v\|_{L^2_y L^\infty_{xT}} + \|v\|_{L^3_T L^\infty_{xy}} \equiv \sum_{i=1}^{10} \eta_i(v).\label{norm}
\end{align}

In the first place, we will consider $s\in(3/4,1]$. Assume that $v_0 \in H^s(\mathbb R^2)$. By using the Strichartz type, local type, and maximal type estimates (inequalities \eqref{Stri1_x} to \eqref{MTE_y}), it can be seen that
\begin{align}
\eta_s(V(t)v_0)\leq c(s)(1+T)^{1/2} \|v_0\|_{H^s_{xy}}.\label{EstLP}
\end{align}

Now, let us assume that $v\in X_{T,s}$. We divide the analysis into the cases $3/4<s\leq 1$, and $s>1$. If $3/4<s<1$, it can be shown that the inequality
\begin{align}
\eta_s(E(v)) \leq c(s)(1+T)^{1/2}(T^{1/6}+T^{2/9}) [\eta_s(v)]^3 \label{EstNLP}
\end{align}
holds, using the Leibniz fractional rule (Lemma \ref{LFR}), Minkowski's inequality, Hölder's inequality, and \eqref{Stri1_x} to \eqref{MTE_y}. The case $s=1$ is simpler since it is used the ordinary Leibniz rule for derivatives. Thus, it is not hard to see that \eqref{EstNLP} is also valid in such case. Consequently, combining \eqref{EstLP} and \eqref{EstNLP}, we have that, for $s\in (3/4,1]$,
\begin{align*}
\eta_s(\Psi(v))\leq c(s)(1+T)^{1/2} \|v_0\|_{H^s_{xy}} + c(s)(1+T)^{1/2}(T^{1/6}+T^{2/9}) [\eta_s(v)]^3.
\end{align*}

Defining $a:=2c(s)(1+T)^{1/2} \|v_0\|_{H^s_{xy}}$, and choosing $T>0$ such that $a^2 c(s)(1+T)^{1/2}(T^{1/6}+T^{2/9})\leq \frac 14$, it follows that
$$\eta_s(v)<a \implies \eta_s(\Psi(v))\leq \frac 12 a + \frac{1}{4a^2} a^3 = \frac 34a < a.$$

If we define $X_{T,s}^a := \{v\in X_{T,s}: \eta_s(v)<a\}$, then we have proven that $\Psi: X_{T,s}^a \to X_{T,s}^a$ is well defined. Moreover, following standard-by-now arguments, it can be seen that $\Psi: X_{T,s}^a \to X_{T,s}^a$ is a contraction. Therefore, there exists a unique $v\in X_{T,s}^a$ such that $\Psi(v)=v$, As a consequence of the method, it can be shown that for every $T'\in(0,T)$, there exists a neighborhood $V$ of $v_0$ in $H^s(\mathbb R^2)$ such that the data-solution map $\tilde v_0\mapsto \tilde v$ from $V$ into $X_{T'}$ is Lipschitz.\\

Now we consider the case $s>1$. For $v_0\in H^s(\mathbb R^2)$, it is evident that $v_0\in H^1(\mathbb R^2)$. Therefore, by applying the previous established result for $s=1$, there exists a unique solution $\tilde v\in X_{T,s}^a$ to the integral equation \eqref{IntEq}. Let us confirm that actually $\tilde v\in C([0,T];H^s(\mathbb R^2))$. Without loss of generality, let us assume that $s\in (1,2]$. Since $\eta_1(\tilde v)<\infty$, it can be seen that
\begin{align}
\|E(\tilde v) \|_{L^\infty_T L^2_{xy}}< cT^{2/9} [\eta_1(\tilde v)]^3<\infty.\label{est_bigS1}
\end{align}

By employing Lemma \ref{Int_Der}, and Hölder's inequality, we can affirm that
\begin{align}
\|D^s_x E(\tilde v)\|_{L^\infty_T L^2_{xy}}\leq c \|\partial_x E(\tilde v)\|_{L^\infty_T L^2_{xy}}^{1-(s-1)/s}\|\partial_x [\partial_x E(\tilde v)]\|_{L^\infty_T L^2_{xy}}^{(s-1)/s}.\label{est_bigS2}
\end{align}

Taking into account the case $s=1$, we notice that
\begin{align}
\|\partial_x E(\tilde v)\|^{(1-(s-1)/2)}_{L^\infty_TL^2_{xy}}<c\{(T^{1/6}+T^{2/9})\eta_1(\tilde v)\}^{3(1-(s-1)/2)}<\infty.\label{est_bigS3}
\end{align}

In relation to the second factor on the right hand side of \eqref{est_bigS2}, we can use estimate \eqref{LTE_xd} in Lemma \ref{LTE}, and the fact that $v\in X_{T,1}^a$ to conclude that
\begin{equation}\label{est_bigS4}
    \begin{split}
        \|\partial_x [\partial_x E(\tilde v)(t)]\|_{L^\infty_TL^2_{xy}}^{(s-1)/s}&\leq c\|\partial_x(v^2 v_x + v^2 v_y) \|_{L^1_x L^2_{yT}}^{(s-1)/s}\\&\hspace{-25mm}\leq c\|v v_x^2\|_{L^1_x L^2_{yT}}^{(s-1)/s}+ c\|v^2 v_{xx}\|_{L^1_x L^2_{yT}}^{(s-1)/s}+ c\| v v_x v_y\|_{L^1_x L^2_{yT}}^{(s-1)/s}+ c\|v^2 v_{xy}\|_{L^1_x L^2_{yT}}^{(s-1)/s}\\&\hspace{-25mm}\leq c[T^{1/18}\|v\|_{L^2_x L^\infty_{yT}} \|v_x\|_{L^\infty_T L^2_{xy}} \|v_x\|_{L^2_T L^\infty_{xy}}]^{(s-1)/s}+ c[\|v\|_{L^2_x L^\infty_{yT}} \|v\|_{L^2_x L^\infty_{yT}} \|v_{xx}\|_{L^\infty_x L^2_{yT}}]^{(s-1)/s}\\&\hspace{-25mm}+ c[T^{1/18}\|v\|_{L^2_xL^\infty_{yT}} \|v_x\|_{L^\infty_T L^2_{xy}} \|v_y\|_{L^2_T L^\infty_{xy}} ]^{(s-1)/s} + [\|v\|_{L^2_x L^\infty_{yT}} \|v\|_{L^2_x L^\infty_{yT}} \|v_{xy}\|_{L^\infty_x L^2_{yT}}]^{(s-1)/s}\\ &\hspace{-25mm}\leq c (1+T^{1/18}) [\eta_1(\tilde v)]^{3(s-1)/s}.
    \end{split}
\end{equation}

From \eqref{est_bigS2}, \eqref{est_bigS3}, and \eqref{est_bigS4} we have that
\begin{align}
\| D^s_x E(\tilde v)\|_{L^\infty_T L^2_{xy}}< c(T^{1/6})^{1-(s-1)/s}(1+T^{1/18}) [\eta_1(\tilde v)]^3<\infty.\label{est_bigS5}
\end{align}

Similary we can obtain that
\begin{align}
\| D^s_y E(\tilde v)\|_{L^\infty_T L^2_{xy}}< c(T^{1/6})^{1-(s-1)/s}(1+T^{1/18}) [\eta_1(\tilde v)]^3<\infty.\label{est_bigS6}
\end{align}

Combining \eqref{est_bigS1}, \eqref{est_bigS5}, and \eqref{est_bigS6}, we conclude that $\tilde u\in C([0,T]; H^s(\mathbb R^2))$. Theorem \ref{LWPHS} is proven.\\

\end{proof}

After establishing the result in Theorem \ref{LWPHS}, we proceed to demonstrate a similar one in weighted Sobolev spaces of the form
$$Z_{s,(r_1,r_2)}:=H^s(\mathbb R^2)\cap L^2((1+|x|^{2r_1}+|y|^{2r_2})dxdy).$$

For this purpose, we will employ a result by Fonseca and Pachón in \cite{FP2016}, adapted for the group associated with the linear part of the symmetrized ZK equation, defined in \eqref{G}, instead of the non-symmetrized equation. It is worth noting that the same proof conducted by Fonseca and Pachón applies to the symmetrized group, with the respective adjustments required. The result is stated in the following Lemma.

\begin{lemma}\label{FP2016} Let $r_1,r_2\in(0,1)$, $s\geq 2\max \{r_1,r_2\}$, and $\{V(t):t\in\mathbb R\}$ the unitary group defined in \eqref{G}. If $u_0\in Z_{s,(r_1,r_2)}$, then, for all $t\in\mathbb R$, and almost every $(x,y)\in\mathbb R^2$,
\begin{align}
|x|^{r_1} V(t) u_0(x,y)= V(t)(|x|^{r_1} u_0)(x,y) + V(t) \{ \Phi_{1,t,r_1}(\hat u_0)(\xi,\eta)\}^{\vee}(x,y),\label{FP2016_1}
\end{align}
with
\begin{align}\label{FP2016_2}
\| \{\Phi_{1,t,r_1}(\hat u_0)(\xi,\eta)\}^{\vee}\|_{L^2_{xy}} \leq c (1+|t|)(\|u_0\|_{L^2_{xy}}+ \|D^s_x u_0\|_{L^2_{xy}}+ \|D^s_y u_0\|_{L^2_{xy}}),
\end{align}
and
\begin{align}
|y|^{r_2} V(t) u_0(x,y)= V(t)(|y|^{r_2} u_0)(x,y) + V(t) \{ \Phi_{2,t,r_2}(\hat u_0)(\xi,\eta)\}^{\vee}(x,y),\label{FP2016_3}
\end{align}
with
\begin{align}
\| \{\Phi_{2,t,r_2}(\hat u_0)(\xi,\eta)\}^{\vee}\|_{L^2_{xy}} \leq c (1+|t|)(\|u_0\|_{L^2_{xy}}+ \|D^s_x u_0\|_{L^2_{xy}}+ \|D^s_y u_0\|_{L^2_{xy}}).\label{FP2016_4}
\end{align}
Moreover, if in addition we suppose that for $\beta\in(0,\min\{r_1,r_2\})$
\begin{equation}
    \label{FP2016_5} D^{\beta}(|x|^{r_1}u_0), D^{\beta}(|x|^{r_1}u_0)\in L^2(\R^2) \ \mbox{ and }\ u_0\in H^{\beta+s}(\R^2),
\end{equation}
then for all $t\in\R$ and for almost every $(x,y)\in\R^2$
\begin{equation}
    \label{FP2016_6} D^\beta(|x|^{r_1}V(t)u_0)(x,y)=V(t)(D^\beta|x|^{r_1}u_0)(x,y)+V(t)D^\beta(\{\Phi_{1,t,r_1}(\widehat{u}_0)(\xi,\eta) \}^\vee)(x,y)
\end{equation}
and 
\begin{equation}
    \label{FP2016_7} D^\beta(|y|^{r_2}V(t)u_0)(x,y)=V(t)(D^\beta|y|^{r_2}u_0)(x,y)+V(t)D^\beta(\{\Phi_{2,t,r_2}(\widehat{u}_0)(\xi,\eta) \}^\vee)(x,y)
\end{equation}
with
\begin{equation}
    \label{FP2016_8}\|D^\beta (\{\Phi_{j,t,r_j}(\widehat{u}_0(\xi,\eta) \}^\vee) \|_{L^2_{xy}}\le c(1+|t|)(\|u_0\|_{L^2_{xy}}+\|D^{\beta+s}_x u_0 \|_{L^2_{xy}}+\|D^{\beta+s}_yu_0 \|_{L^2_{xy}}), \ \ j=1,2.
\end{equation}

\end{lemma}
\begin{proof}
    See Theorem 2 in \cite{FP2016}.
\end{proof}
The result in weighted Sobolev spaces is as follows.

\begin{theorem}\label{LWPZ}
Let $s>3/4$ and $v\in H^s(\mathbb R^2)$ the solution of the IVP \eqref{IVP_smZK} obtained in Theorem \ref{LWPHS}. Let us assume that $ (|x|^{r_1}+|y|^{r_2})v_0\in L^2(\mathbb R^2)$, with $0<2\max\{r_1,r_2\}\leq s$. Then $v\in C([0,T];Z_{s,(r_1,r_2)})$.\\

Moreover, for any $T'\in(0,T)$ there exists a neighborhood $V$ of $v_0$ in $Z_{s,(r_1,r_2)}$ such that the data-solution map $\tilde v_0 \mapsto \tilde v$ from $V$ into the class defined by $X_{T,s}$ in Theorem \ref{LWPHS}, with $T'$ instead of $T$, is Lipschitz.
\end{theorem}

\begin{proof} We follow the ideas from the proof of Theorem 2 in \cite{FP2016}. We begin, as in the proof of Theorem \ref{LWPHS}, assuming that $s\in (3/4,1]$. In the proof of Theorem \ref{LWPHS} we obtained a unique $\tilde v\in X_{T,s}\subset C([0,T];H^s(\mathbb R^2))$ such that
\begin{align*}
\tilde v(t):= V(t) v_0 - 4^{-1/3} \int_{-t}^t V(t-t')(\tilde v \tilde v_x + \tilde v \tilde v_y)(t') dt'\equiv V(t)v_0 + E(\tilde v)(t),
\end{align*}
and we proved that 
\begin{align}
\eta_s(E(v)) \leq & c(s)(1+T)^{1/2}( T^{1/6}+T^{2/9}) [\eta_s(v)]^3,\label{est_norm1} \\ 
\eta_s(\Psi(v))\leq & c(s)(1+T)^{1/2}\|v_0\|_{H^s(x,y)} + c(s)(1+T)^{1/2}( T^{1/6}+T^{2/9} )[\eta_s(v)]^3,\label{est_norm2}
\end{align}
for every $v\in X_{T,s}$.\\

Let us define the norm
$$\eta_{r_1,r_2}(v):= \| |x|^{r_1} v\|_{L^\infty_T L^2_{xy}} + \| |y|^{r_2} v\|_{L^\infty_T L^2_xy},$$
By using Lemma \ref{FP2016}, Hölder's inequality, and bearing in mind \eqref{est_norm1}, we can affirm that
\begin{align*}
\| |x|^{r_1} \tilde v(t) \|_{L^2_{xy}} \leq &\| |x|^{r_1} V(t)v_0\|_{L^2_{xy}} + c \left\| |x|^{r_1} E(\tilde v)(t) \right\|_{L^2_{xy}}\\
\leq  &\|V(t) (|x|^{r_1} v_0)\|_{L^2_{xy}} + c(1+T) \|v_0\|_{H^s_{xy}} + c \int_0^T \| |x|^{r_1} (\tilde v^2\tilde v_x)(t') \|_{L^2_{xy}} dt'\\
& + c \int _0^T \| |x|^{r_1} (\tilde v^2\tilde v_y)(t') \|_{L^2_{xy}} dt' + c (1+T) \| E(\tilde v)(t) \|_{H^s_{xy}}\\
\leq & \| |x|^{r_1} v_0\|_{L^2_{xy}} + c(1+T) \| v_0\|_{H^s_{xy}} + c T^{2/9} \| |x|^{r_1} \tilde v\|_{L^\infty_T L^2_{xy}} \|\tilde v\|_{L^3_T L^\infty_{xy}} \|\tilde v_x\|_{L^{9/4}_T L^\infty_{xy}}\\
& + c T^{2/9} \| |x|^{r_1} \tilde v\|_{L^\infty_T L^2_{xy}} \|\tilde v\|_{L^3_T L^\infty_{xy}} \|\tilde v_y\|_{L^{9/4}_T L^\infty_{xy}} + c(1+T)(T^{1/6}+T^{2/9})[\eta_s(\tilde v)]^3,
\end{align*}

which implies that 
\begin{equation*}
    \begin{split}
        \| |x|^{r_1} \tilde v(t))\|_{L^\infty_T L^2_{xy}}& \leq \| |x|^{r_1} v_0\|_{L^2_{xy}} + c(1+T) \|v_0\|_{H^s_{xy}} + cT^{2/9} \| |y|^{r_2} \tilde v \|_{L^\infty_T L^2_{xy}} [\eta_s(\tilde v)]^2 \\&\hspace{5mm}+ c(1+T)(T^{1/6}+T^{2/9}) [\eta_s(\tilde v)]^3.
    \end{split}
\end{equation*}

In a similar way,
\begin{equation*}
    \begin{split}
        \| |y|^{r_2} \tilde v(t))\|_{L^\infty_T L^2_{xy}}& \leq \| |y|^{r_2} v_0\|_{L^2_{xy}} + c(1+T) \|v_0\|_{H^s_{xy}} + cT^{2/9} \| |y|^{r_2} \tilde v \|_{L^\infty_T L^2_{xy}} [\eta_s(\tilde v)]^2 \\&\hspace{5mm}+ c(1+T)(T^{1/6}+T^{2/9}) [\eta_s(\tilde v)]^3.
    \end{split}
\end{equation*}

Thus
\begin{equation*}
    \begin{split}
\eta_{r_1,r_2}(\tilde v)(t)&\leq \| |x|^{r_1} v_0\|_{L^2_{xy}} + \| |y|^{r_2} v_0\|_{L^2_{xy}} + c(1+T) \|v_0\|_{H^s_{xy}} + cT^{2/9}  \eta(\tilde v) [\eta_s(\tilde v)]^2 \\&\hspace{5mm}+ c(1+T)(T^{1/6}+T^{2/9}) [\eta_s(\tilde v)]^3.
    \end{split}
\end{equation*}

If we define
\begin{align*}
\tilde\eta_{s,(r_1,r_2)}(v):=\eta_s(v)+\eta_{r_1,r_2}(v),\label{norma_peso}
\end{align*}
we have that
\begin{equation*}
    \begin{split}
        \tilde \eta_{s,(r_1,r_2)}(\tilde v)&=\eta_s(\tilde v)+\eta_{r_1,r_2}(\tilde v)\leq  \eta_s(\tilde v)+c[\| |x|^{r_1} v_0\|_{L^2_{xy}} + \| |y|^{r_2} v_0\|_{L^2_{xy}} +  (1+T)\|v_0\|_{H^s_{xy}}] \\&\hspace{5mm}+ cT^{2/9}\tilde \eta_{s,(r_1,r_2)}(\tilde v)[\eta_s(\tilde v)]^2 + c(1+T)(T^{1/6}+T^{2/9}) [\eta_s(\tilde v)]^3\\&\le \eta_s(\tilde{v})+c[\||x|^{r_1}v_0\|_{L^2_{xy}}+\||y|^{r_2}v_0\|_{L^2_{xy}}+(1+T)\|v_0\|_{H^s_{xy}} ]\\&\hspace{5mm}+c(1+T)^{1/2}(T^{1/6}+T^{2/9})\tilde\eta_{s,(r_1,r_2)}(\tilde{v})[\eta_s(\tilde{v})]^2+c(1+T)(T^{1/6}+T^{2/9})[\eta_s(\tilde{v})]^3.
    \end{split}
\end{equation*}

Let us recall that in the proof of Theorem \ref{LWPHS}, we chose $a:=2c(s)(1+T)^{1/2} \|v_0\|_{H^s_{xy}}$ and  $T>0$ such that $a^2 c(s)(1+T)^{1/2}(T^{1/6}+T^{2/9})\leq \frac 14$. Besides, from the contraction principle, we had that $\eta_s(\tilde v)<a$. Therefore
\begin{equation*}
    \tilde \eta_{s,(r_1,r_2)}(\tilde v) \leq  a+c[\| |x|^{r_1} v_0\|_{L^2_{xy}} + \| |y|^{r_2} v_0\|_{L^2_{xy}} +  (1+T)\|v_0\|_{H^s_{xy}}] + \frac14 \tilde \eta_{s,(r_1,r_2)}(\tilde v) + c(1+T)^{1/2}a.
\end{equation*}
The latter implies
\begin{equation*}
    \tilde \eta_{s,(r_1,r_2)}(\tilde v)\leq c[\| |x|^{r_1} v_0\|_{L^2_{xy}} + \| |y|^{r_2} v_0\|_{L^2_{xy}} +  (1+T)\|v_0\|_{H^s_{xy}}]+c(1+T)^{1/2}a<\infty.
\end{equation*}

This completes the proof for $s\in(3/4,1]$. The case $s>1$ can be handled in a similar manner to that of the proof of Theorem \ref{LWPHS} taking into account that \eqref{FP2016_2} and \eqref{FP2016_4} can be adapted for $r_1,r_2\ge 1$ via iterations at every integer step (as done in the proof of Lemma \ref{remark2} below). Theorem \ref{LWPZ} is proven.

\end{proof}

\begin{lemma}\label{remark1}
    Let $s>9/4$ and $r_1,r_2>9/8$ such that $2\max\{r_1,r_2\}\leq s$. If $v\in C([0,T];Z_{s,(r_1,r_2)})$ is the solution of the IVP \eqref{IVP_smZK} provided by Theorem \ref{LWPZ}, then \begin{equation}
        \|\langle x\rangle^{\frac{3}{4}^+}v\|_{L^4_xL^\infty_{yt}}+\|\langle y\rangle^{\frac{3}{4}^+}v\|_{L^4_yL^\infty_{xt}}\le c(T,\eta_{9/4+}(v),\|v\|_{L^\infty_TZ_{s,(r_1,r_2)}}).
    \end{equation}
    (See the definition of $\eta_{9/4+}(v)$ in \eqref{EstNLP}).
\end{lemma}
\begin{proof}
    Using Lemma \ref{FP2016} and \eqref{grunrockcol} we have
    \begin{equation}
        \label{D1}
        \begin{split}
            \|\langle x\rangle^{\frac{3}{4}^+}V(t)v_0\|_{L^4_xL^\infty_{yt}}&\le \|V(t)\langle x \rangle^{\frac{3}{4}^+}v_0\|_{L^4_xL^\infty_{yt}}+\|V(t)\{\Phi_{1,t,\frac{3}{4}^+}(\widehat{v}_0)\}^\vee\|_{L^4_xL^\infty_{yt}}\\&\le c\|J_{xy}^{\frac{3}{4}^+}\langle x\rangle^{\frac{3}{4}^+}v_0\|_{L^2_{xy}}+c\|J_{xy}^{\frac{3}{4}^+}\{\Phi_{1,t,\frac{3}{4}^+}(\widehat{v}_0 )\}^\vee\|_{L^2_{xy}}.
        \end{split}
    \end{equation}
Note from \eqref{interpol1} with $\theta=1/3$ that
\begin{equation}
    \label{D2} \|J^{\frac{3}{4}^+}_{xy}\langle x\rangle^{\frac{3}{4}^+}v_0 \|_{L^2_{xy}}\le c\|J^{\frac{9}{4}^+}_{xy}v_0\|_{L^2_{xy}}^{1/3}\|\langle x\rangle^{\frac{9}{8}^+}v_0\|_{L^2_{xy}}^{2/3} \le c(\|v_0\|_{Z_{s,(r_1,r_2)}}).
\end{equation}
Similarly, from \eqref{FP2016_2} and \eqref{FP2016_8} we have
\begin{equation}
    \label{D3}\|J_{xy}^{\frac{3}{4}^+}\{\Phi_{1,t,\frac{3}{4}^+}(\widehat{v}_0) \}^\vee\|_{L^2_{xy}}\le c(1+T)\|J^{\frac{9}{4}^+}_{xy}v_0\|_{L^2_{xy}}\le c(T, \|v_0\|_{Z_{s,(r_1,r_2)}}).
\end{equation}
From \eqref{D1} - \eqref{D3} we conclude
\begin{equation}
    \label{D4}\|\langle x\rangle^{\frac{3}{4}^+}V(t)v_0\|_{L^4_xL^\infty_{yt}}\le c(T, \|v_0\|_{Z_{s,(r_1,r_2)}}).
\end{equation}

On the other hand, 
\begin{equation}\label{D5}
    \begin{split}
        \int_0^T\|\langle x\rangle^{\frac{3}{4}^+}V(t-t')(v^2(\partial_x v+\partial_y v)&)(t')\|_{L^4_xL^\infty_{yt}}dt'\\&\hspace{-15mm}\le \int_0^T\|V(t-t')\langle x\rangle^{\frac{3}{4}^+}(v^2(\partial_x v+\partial_y v))(t')  \|_{L^4_xL^\infty_{yt}}dt'\\&\hspace{-10mm}+\int_0^T \|V(t-t')\{\Phi_{1,t-t',\frac{3}{4}^+}(v^2(\partial_x v+\partial_y v))^\wedge \}^\vee \|_{L^4_xL^\infty_{yt}}dt'\\&\hspace{-15mm}\le c\int_0^T \|J^{\frac{3}{4}^+}_{xy}\langle x\rangle^{\frac{3}{4}^+}(v^2(\partial_x v\partial_y v))(t') \|_{L^2_{xy}}dt'\\&\hspace{-10mm}+c\int_0^T\|J^{\frac{3}{4}^+}_{xy}\{\Phi_{1,t-t,\frac{3}{4}^+}(v^2(\partial_x v+\partial_y v))^\wedge \}^\vee \|_{L^2_{xy}}dt'.
    \end{split}
\end{equation}
Note as done for \eqref{D2}, via interpolation, and using Young's inequality, we have that
\begin{equation}
    \label{D6} 
    \begin{split}
        \int_0^T\|J^{\frac{3}{4}^+}_{xy}\langle x\rangle^{\frac{3}{4}^+} (v^2(\partial_xv+\partial_y v))(t')\|_{L^2_{xy}}dt'&\le c\int_0^T\|J^{\frac{9}{4}^+}_{xy}(v^2(\partial_x v+\partial_y v))(t')\|_{L^2_{xy}}dt'\\&\hspace{3mm}+c\int_0^T\|\langle x\rangle^{\frac{9}{8}^+}(v^2(\partial_x v+\partial_y v))(t')\|_{L^2_{xy}}dt'\\&\le c(T,\eta_{9/4+}(v),\|v\|_{L^\infty_TZ_{s,(r_1,r_2)}}).
    \end{split}
\end{equation}
 Also
 \begin{equation}\label{D7}
 \begin{split}
     \int_0^T \|J^{\frac{3}{4}^+}_{xy}\{\Phi_{1,t-t',\frac{3}{4}^+}(v^2(\partial_x v+\partial_y v)) \} \|_{L^2_{xy}}dt'&\le c(1+T)\int_0^T\|v^2(\partial_x v+\partial_y v)(t')\|_{H^{\frac{9}{4}^+}_{xy}}dt'\\&\le c(T,\eta_{9/4+(v)}). 
 \end{split}
 \end{equation}

From \eqref{D5}-\eqref{D7} we conclude \begin{equation}
    \label{D8} \int_0^T \|\langle x\rangle^{\frac{3}{4}^+}V(t-t')(v^2(\partial_x v+\partial_y v))(t') \|_{L^4_xL^\infty_{yt}}dt'\le c(T,\vertiii{v},\|v\|_{L^\infty_TZ_{s,(r_1,r_2)}}).
\end{equation}

Finally, from \eqref{D4}, \eqref{D8} and symmetry of the arguments we obtain
\begin{equation*}
    \|\langle x\rangle^{\frac{3}{4}^+}v \|_{L^4_xL^\infty_{yt}}+\|\langle y\rangle^{\frac{3}{4}^+}v \|_{L^4_yL^\infty_{xt}}\le c(T,\eta_{9/4+}(v),\|v\|_{L^\infty_TZ_{s,(r_1,r_2)}},\|v_0\|_{Z_{s,(r_1,r_2)}}),
\end{equation*}

which proves Lemma \ref{remark1}

\end{proof}

\begin{lemma}
    \label{remark2} Let $s>9/4$ and $r_1,r_2>9/8$ such that $2\max\{r_1,r_2\}\leq s$. If $v\in C([0,T];Z_{s,(r_1,r_2)})$ is the solution of the IVP \eqref{IVP_smZK} provided by Theorem \ref{LWPZ}, then
    \begin{equation*}
        \|\langle x\rangle^{s/2} v\|_{L^3_TL^\infty_{xy}}+ \|\langle y\rangle^{s/2} v\|_{L^3_TL^\infty_{xy}}\le c(T,\eta_{s}(v),\|v\|_{L^\infty_TZ_{s,s/2}},\|v_0\|_{Z_{s,(r_1,r_2)}}).
    \end{equation*}
    (See the definition of $\eta_{s}(v)$ in \eqref{EstNLP}).
\end{lemma}
\begin{proof}
It is sufficient to develop the proof under the assumption $\frac{s}{2}=1+b$ with $b\in(0,1)$. 

Given that $\Gamma_{x,t}:=x-3t\partial_x^2$ commutes with $\partial_t+\partial_x^3+\partial_y^3$, it can be seen that \begin{equation*}
    xV(t)(v_0)(x,y)=V(t)(xv_0)(x,y) + 3tV(t)(\partial_x^2v_0)(x,y).
\end{equation*}
Thus, \begin{equation}
    \begin{split}
        |x|^{s/2} |V(t)(v_0)(x,y)|&\leq |x|^b |V(t)(xv_0)(x,y)+3tV(t)(\partial_x^2v_0)(x,y)|\\& \leq |V(t)(|x|^b x v_0)(x,y)| + |V(t)\{\Phi_{1,t,b}(\widehat{xv_0}) \}^\vee(x,y)| \\&\hspace{5mm} + |3tV(t)(|x|^b\partial_x^2 v_0)(x,y)| + |3tV(t)\{\Phi_{1,t,b}(\widehat{\partial_x^2v_0}) \}^\vee(x,y)|.\label{ineq_3.29}
    \end{split}
\end{equation}

Therefore, taking into account the Strichartz estimate \eqref{Stri1_x} with $\epsilon=0$, inequality \eqref{ineq_3.29}, estimate \eqref{FP2016_2}, and the interpolation inequalities \eqref{interpol1r} and \eqref{interpol1} with $\theta=\frac{2b}s$, we have

\begin{equation}\label{p1}
    \begin{split}
        \|\langle x\rangle^{s/2}V(t)v_0\|_{L^3_TL^\infty_{xy}}&\le c(T) [ \|\langle x\rangle^{s/2}v_0\|_{L^2_{xy}}+\|\{\Phi_{1,t,b}(\widehat{xv_0})\}^\vee \|_{L^2_{xy}} + \|\langle x\rangle^{b}\partial_x^{2}v_0\|_{L^2_{xy}}\\&\hspace{5mm}+ \|\{\Phi_{1,t,b}(\widehat{\partial_x^2 v})\}^{\vee} \|_{L^2_{xy}}] \\&\le c(T) [\|\langle x\rangle^{s/2}v_0\|_{L^2_{xy}}+ \|J^{2b}_{xy}\langle x\rangle v_0 \|_{L^2_{xy}}+ \|\langle x\rangle^b J_x^2v_0\|_{L^2_{xy}}\\&\hspace{5mm} + \|J_{xy}^{2b}\partial_x^2v_0\|_{L^2_{xy}}] \\&\le c(T) [\|\langle x\rangle^{s/2}v_0\|_{L^2_{xy}}+ \|J^{s}_{xy}v_0\|_{L^2_{xy}}^{\theta}\|\langle x\rangle^{s/2}v_0\|_{L^2_{xy}}^{1-\theta}\\&\hspace{5mm}+ \|J^{s}_{xy}v_0\|_{L^2_{xy}}^{1-\theta}\|\langle x\rangle^{s/2}v_0\|_{L^2_{xy}}^{\theta}+ \|J^s_{xy}v_0\|_{L^2_{xy}}]\\&\le c(T)\|v_0\|_{Z_{s,(r_1,r_2)}}.
    \end{split}
\end{equation}

A similar procedure shows that the nonlinear part satisfies
\begin{equation}\label{p2}
    \begin{split}
        \left\|\langle x\rangle^{s/2}\int_0^tV(t-t')(v^2(\partial_x v+\partial_y v))(t')dt'\right\|_{L^3_TL^\infty_{xy}}&\\& \hspace{-60mm}\le c(T)\int_0^T\left\{\|\langle x\rangle^{s/2}(v^2(\partial_x v+\partial_y v))(t')\|_{L^2_{xy}}+ \|(v^2(\partial_x v+\partial_y v))(t')\|_{H^s_{xy}} \right\}dt'\\&\hspace{-60mm}\le c(T,\eta_s(v),\|v\|_{ Z_{s,s/2}}).
    \end{split}
\end{equation}

The lemma follows from \eqref{p1}, \eqref{p2} and symmetry of the arguments.

\end{proof}


\section{Decay versus regularity}
For $N\in \N$ and $\alpha\in[0,1/2)$ denote 
\begin{equation*}
    \Tilde{\Tilde{\rho}}_{N}^\alpha (x):=
    \begin{cases}
        \langle x \rangle^{2\alpha}-1 \hspace{5mm}  x\in [0,N],\\
        (2N)^{2\alpha} \hspace{9mm}  x\in [3N,\infty).
    \end{cases}
\end{equation*}
Let $\Tilde{\rho}_{N}^\alpha$ be a smooth regularization of $\Tilde{\Tilde{\rho}}_{N}^\alpha$, defined in $[0,\infty)$, such that for $j=1,2,\dots$ we have $\left|\partial_x^j \Tilde{\rho}_{N}^\alpha\right|\le c$ with $c$ independent of $N$. Set $\rho_{N}^\alpha$ to be the odd extension of $\Tilde{\rho}_{N}^\alpha$, that is, \begin{equation}
    \label{c0}
    \rho_N^\alpha(x):=\begin{cases}
        \Tilde{\rho}_{N}^\alpha(x) \hspace{11mm}  x\in [0,\infty],\\
        -\Tilde{\rho}_{N}^\alpha(-x) \hspace{5mm}  x\in (-\infty,0).
    \end{cases}
\end{equation}
Note that $(\rho_N^\alpha)^\prime\ge 0$.\\

Let $m\in \N\cup\{0\}$. For $\beta>m/2$ define $\rho_{N,m}^\beta:=\rho_N^{\beta-m/2}$. Let us note that in case $\beta\in(\frac{m+k}{2},\frac{m+k+1}{2}]$ for some $k\in \N$ we would have that $\beta-m/2\in (\frac{k}{2},\frac{k+1}{2}]$ which prevents the derivatives $|\partial_x^j \rho_{N,m}^\beta|$, $j=1,\dots,k$; to be bounded above independent on $N$.

Before stating and proving the main results of this work, let us recall the definition of the norm $\eta_s$ given in \eqref{norm}, as it will be used frequently in the proofs that follow.

\begin{align*}
\notag \eta_s(v):= & \|v\|_{L^\infty_T H^s_{xy}} + \| D^s_x v_x\|_{L^\infty_x L^2_{yT}} + \|D^s_y v_x\|_{L^\infty_x L^2_{yT}} + \|v_x\|_{L^{9/4}_T L^\infty_{xy}} + \|v\|_{L^2_x L^\infty_{yT}} + \|D^s_x v_y\|_{L^\infty_y L^2_{xT}}\\
& + \|D^s_y v_x\|_{L^\infty_y L^2_{xT}} + \|v_y\|_{L^{9/4}_T L^\infty_{xy}} + \|v\|_{L^2_y L^\infty_{xT}} + \|v\|_{L^3_T L^\infty_{xy}} \equiv \sum_{i=1}^{10} \eta_i(v),
\end{align*}

\begin{theorem}
    \label{teo1}
    Let $v\in C([0,T];H^{\frac{3}{4}^+}(\R^2))$ be the solution of the IVP \eqref{IVP_smZK} provided by Theorem \ref{LWPHS}. Let us assume that there are $t_0<t_1 \in [0,T]$ such that 
    \begin{equation}
        (1+|x|^{\frac{3}{8}^++\alpha_1}+|y|^{\frac{3}{8}^++\alpha_2})v(t_i)\in L^2(\R^2), \ i=0,1,    \end{equation}
    where $\alpha_1,\alpha_2\in [0,1/8)$ and denote $\alpha:=\min\{\alpha_1,\alpha_2\}$. Then, $v\in C([0,T];H^{\frac{3}{4}^++2\alpha}(\R^2))$.
\end{theorem}

\begin{proof}
Without loss of generality we assume $t_0=0$.  Let $\{v_{0m}\}_m\subset C_0^\infty(\R)$ be a sequence converging to $v_0$ in $H^{\frac{3}{4}^+}(\R^2)$. Denote with $v_m \in H^\infty(\R^2)$ the corresponding solution provided by the local theory with initial datum $v_{0m}$. By regularity of the data-solution map we can assume all the $v_m$'s are defined in $[0,T]$ with $v_m$ converging to $v$ in $C\left([0,T];H^{\frac{3}{4}^+}(\R^2)\right)$. 

We set $\beta_1:=\frac{3}{8}^+ +\alpha_1$ and consider $\rho_N\equiv \rho_{N,0}^{\beta_1}$ defined above. In particular we note that $|\rho_N(x)|\le c\langle x\rangle^{2\beta_1}$. 

Let us take the symmetrized mZK equation \eqref{smZK} for $v_m$ and multiply it by $v_m(x,y) \rho_N(x)$. Taking into account that $\rho_N$ depends only on $x$, after integration by parts in space we get

\begin{equation}\label{a1}
    \begin{split}
        \frac{1}{2}\frac{d}{dt}\int v_m^2 \rho_N dxdy +\frac{3}{2}\int (\partial_x v_m)^2\rho_N^\prime dxdy - \frac{1}{2}\int v_m^2\rho_N^{\prime \prime \prime}dxdy -\frac{\mu}{4}\int v_m^4 \rho_N'dxdy=0.
    \end{split}
\end{equation}

In what follows, we will make use several times (without explicitly mentioning it) the continuous dependence provided by Theorem \ref{LWPHS} in the sense that for $m$ large enough $\eta_s(v_m)\leq \eta_s(v)$. 

Integrating \eqref{a1} over $[0,t_1]$ and taking into account that 
\begin{equation}\label{a2}
    \begin{split}
        \left|\int_0^{t_1}\int v_m^2 \rho_N^{\prime\prime\prime}dxdydt \right|&\le c{t_1}\|v_m\|^2_{L^\infty_TL^2_{xy}}\le c\|v\|^{2}_{L^\infty_TL^2_{xy}}\\&=c(\eta_{3/4+}(v)),
    \end{split}
\end{equation}
and that 
\begin{equation}\label{a3}
    \begin{split}
        \left|\int_0^{t_1}\int v_m^4 \rho_N^\prime dxdydt \right|&\le c\int_0^{t_1}\|v_m\|^{2}_{L^\infty_{xy}}\int v_m^2 dxdydt \\&\le c\|v_m\|_{L^2_TL^\infty_{xy}}^2\|v_m\|_{L^\infty_{T}L^2_{xy}}\\&\le c(T)\|v_m\|^2_{L^3_TL^\infty_{xy}}\|v_m\|_{L^\infty_TL^2_{xy}}\\&\le c(T) [\eta_{3/4+}(v)]^4\le c(T,\eta_{3/4+}(v)),
    \end{split}
\end{equation}
we obtain for $m$ large enough that
\begin{equation}
    \label{a4}
    \begin{split}
        \int_0^{t_1}\int (\partial_xv_m)^2\rho_N^\prime dxdydt &\le c(T,\eta_{3/4+}(v))+c\|v_m(t_1)\rho_N^{1/2}\|^2_{L^2_{xy}}+c\|v_{0m}\rho_N^{1/2}\|^2_{L^2_{xy}}\\&\le 
        c(T,\eta_{3/4+}(v))+c\|v(t_1)\rho_N^{1/2}\|^2_{L^2_{xy}}+c\|v_{0}\rho_N^{1/2}\|^2_{L^2_{xy}}\\&\le 
        c(T,\eta_{3/4+}(v))+c\|v(t_1)\langle x\rangle^{\beta_1}\|^2_{L^2_{xy}}+c\|v_{0}\langle x\rangle^{\beta_1}\|^2_{L^2_{xy}}\\&\le c(T,\eta_{3/4+}(v),\|v(t_1)\|_{L^2_{w}},\|v_0\|_{L^2_w}),
    \end{split}
\end{equation}
where $\|f\|_{L^2_w}:=\|(1+|x|^{\frac{3}{8}^++\alpha_1}+|y|^{\frac{3}{8}^++\alpha_2})f\|_{L^2_{xy}}$. 

Let us note that for each $N$ fixed we have
\begin{equation}
    \label{a5}
    \begin{split}
        \int_0^{t_1}\int \left|(\partial_x v_m)^2 \rho_N'-(\partial_x v)^2\rho_N' \right|dxdydt & \le \int_0^{t_1}\int |\partial_x v_m - \partial_x v||\partial_x v_m+\partial_x v|\rho_N^\prime(x)dxdydt\\ &\hspace{-25mm}\le c\int_0^{t_1} \left\|(\partial_x v_m-\partial_x v) \chi_{[-2N,2N]}(\cdot_x) \right\|_{L^2_{xy}}\left\|(\partial_x v_m + \partial_x v)\chi_{[-2N,2N]}(\cdot_x) \right\|_{L^2_{xy}}dt\\ &\hspace{-25mm}\le c\left\|(\partial_x v_m -\partial_x v)\chi_{[-2N,2N]}(\cdot_x) \right\|_{L^2_{xyt_1}}\left\| (\partial_x v_m + \partial_x v)\chi_{[-2N,2N]}(\cdot_x)\right\|_{L^2_{xyt_1}}\\&\hspace{-25mm}\le c(N)\|\partial_x v_m - \partial_x v\|_{L^{\infty}_{x}L^2_{yt_1}}\left(\|\partial_x v_m\|_{L^{\infty}_{x}L^2_{yt_1}}+\|\partial_x v\|_{L^{\infty}_{x}L^2_{yt_1}}\right)\\&\hspace{-25mm}\le c(N,\vertiii{v})\|\partial_x v_m - \partial_x v\|_{L^{\infty}_{x}L^2_{yt_1}}\xrightarrow[m\to\infty]{}0.
    \end{split}
\end{equation}

Therefore, for all $N\in \N$ we have from \eqref{a4} and \eqref{a5} that there exists $c>0$ such that
\begin{equation}\label{a6}
    \int_0^{t_1}\int (\partial_x v)^2\rho_N^\prime dxdydt \le c,
\end{equation}
with $c$ independent of $N$. 

According to the definition of $\rho_N(x)$, for $|x|>1$ it follows that $\rho_N^\prime(x)\xrightarrow[N\to\infty]{}2\beta_1\langle x\rangle^{2\beta_1-2}x \sim \langle x\rangle^{2\beta_1 -1}$. Thus, from Fatou's Lemma, and considering \eqref{a6}, we obtain that
\begin{equation}
    \label{a7}
   \int_0^{t_1}\int \int_{|x|>1} (\partial_x v)^2\langle x\rangle^{2\beta_1 -1}dxdyt\le c \liminf\limits_{N\to\infty}{\int_0^{t_1}\int\int_{|x|>1}(\partial_x v)^2\rho_N^\prime dxdydt }\le c.
 \end{equation}

Let us also note that \begin{equation}
    \label{a8}\begin{split}
    \int_0^{t_1}\int\int_{|x|\le 1}(\partial_x v)^2\langle x\rangle^{2\beta_1 -1}dxdyt& \le c\int_0^{t_1}\int (\partial_x v)^2 \chi_{\{|x|\le 1\}}(\cdot_x)dxdydt\\& \le c \|\partial_x v\|_{L^\infty_xL^2_{yt_1}}\le c(\eta_{3/4+}(v)).
    \end{split}
\end{equation}

Combining \eqref{a7} and \eqref{a8} we conclude that
\begin{equation}
    \label{a9}
    \int_0^{t_1}\int (\partial_x v)^2 \langle x\rangle^{2\beta_1-1}dxdydt <\infty,
\end{equation}
which implies \begin{equation}
    \label{a10}
    \langle x \rangle^{\beta_1-1/2}\partial_x v(t)\in L^2(\R^2) \ a. e.\ t\in [0,t_1].
\end{equation}

Now, with \eqref{a9} in hand, we reapply the argument above with $\phi_N$ being the even extension of $\tilde{\rho}_N$ instead of $\rho_N$. Note that $|\phi_N^\prime|\le c\langle x\rangle^{2\beta_1-1}$. Using \eqref{a1} with $\phi_N$ instead of $\rho_N$ and arguing as in \eqref{a2}-\eqref{a9} it follows that 
\begin{equation*}
    \lim\limits_{m\to\infty}\int_0^{t_1}\int (\partial_x v_m)^2\phi_N^\prime dxdydt =\int_0^{t_1}\int(\partial_x v)^2\phi_N^\prime dxdydt 
\end{equation*}
and 
\begin{equation*}
    \left|\int_0^{t_1}\int (\partial_x v)^2\phi_N^\prime dxdydt \right| \le \int_0^{t_1}\int (\partial_x v)^2\langle x\rangle^{2\beta_1-1}dxdydt \le c.
\end{equation*}

From these estimates and integrating over $[0,t]\subset [0,t_1]$ it follows that
\begin{equation}
    \label{a11} \langle x \rangle^{\beta_1}v(t) \in L^2(\R^2)\ \mbox{for all }t\in[0,t_1]. 
\end{equation}

Because of the symmetry of the equation it can be seen that for $\beta_2:=\frac{3}{8}^++\alpha_2$
\begin{equation}\label{a12}
    \begin{cases}
    \langle y \rangle^{\beta_2} v(t)\in L^2(\R^2)\ \mbox{for all }t\in[0,t_1],\\
    \langle y \rangle^{\beta_2 -1/2}\partial_y v(t) \in L^2(\R^2)\ a. e.\ t\in[0,t_1].
\end{cases}
\end{equation}

Let $t^*\in[0,t_1]$ be such that \eqref{a10}, \eqref{a11} and \eqref{a12} hold at $t^*$. In particular we have \begin{equation*}
    J^1_x(\langle x\rangle^{\beta_1-1/2}v(t^*))\ \mbox{and } J^1_y(\langle y\rangle^{\beta_2-1/2}v(t^*))\ \mbox{are in } L^2(\R^2).
\end{equation*}
For $i=1,2$ set $f_i(x,y):=\langle a_i \rangle^{\beta_i-1/2}v(t^*)(x,y)$ where $a_1=x$ and $a_2=y$. Interpolating with $\theta=2\beta_i\in(0,1)$ we get
\begin{equation}
    \label{a13}
    \left\|J^{2\beta_i}_{a_i}v(t^*) \right\|_{L^2_{xy}}=\left\|J^{\theta \cdot 1}_{a_i}\langle a_i \rangle^{(1-\theta)\cdot\frac{1}{2}}f_i \right\|_{L^2_{xy}}\le c\left\|J^1_{a_i}f_i\right\|_{L^2_{xy}}^{\theta}\left\|\langle a_i\rangle^{1/2}f_i\right\|_{L^2_{xy}}^{1-\theta}<\infty.
\end{equation}

We conclude from \eqref{a13} that $v(t^*)\in H^{\frac{3}{4}^++2\alpha}(\R^2)$, which implies $v\in C([0,T];H^{\frac{3}{4}^++2\alpha}(\R^2))$. Moreover, by Theorem \ref{LWPZ}, \begin{equation}
    v\in C([0,T];Z_{\frac{3}{4}^++2\alpha,\frac{3}{8}^++\alpha}).
\end{equation}

Theorem \ref{teo1} is proven.
\end{proof}

\begin{theorem}\label{T4.2}
    Let $v\in C([0,T];H^{\frac{9}{4}^+}(\R^2))$ be the solution of the IVP \eqref{IVP_smZK} provided by Theorem \ref{LWPHS}. Assume there are $t_0<t_1\in[0,T]$ such that 
    \begin{equation}
        (1+|x|^{\frac{9}{8}^++\alpha_1}+|y|^{\frac{9}{8}^++\alpha_2})v(t_i)\in L^2(\R^2), \ i=0,1;
    \end{equation}
    where $\alpha_1,\alpha_2\ge 0$.  Denote $\alpha:=\min\{\alpha_1,\alpha_2\}$. Then $$v\in C([0,T];H^{\frac{9}{4}^++2\alpha}(\R^2)).$$
\end{theorem}
\begin{proof}
Without loss of generality assume $t_0=0$. We proceed by cases depending on the size of $\alpha$. \\

\noindent \underline{\textbf{Case $\alpha\in [0,3/8)$}.}\\

From the hypothesis and Theorem \ref{LWPZ} it follows that $v\in C([0,T];Z_{\frac{9}{4}^+,\frac{9}{8}^+})$. Let $\{v_{0m}\}_{m}$ be a sequence of $C_0^\infty(\R^2)$ functions converging to $v_0$ in $Z_{\frac{9}{4}^+,\frac{9}{8}^+}$. Let us denote with $v_m$ the respective solution provided by Theorem \ref{LWPZ} with initial datum $v_{0m}$. We can assume that for $m$ large enough all the $v_m$'s are in $C([0,T];Z_{\frac{9}{4}^+,\frac{9}{8}^+})$.\\

In what follows we will make use several times, without mentioning the continuous dependence provided by Theorem \ref{LWPZ} in the sense that for $m$ large enough $\eta_{9/4+}(v_m)\le c \eta_{9/4+}(v)$.\\

We first prove that $\langle x\rangle^{1/2}\partial_x v(t)\in L^2(\R^2)$ \textit{a.e.} in $[0,t_1]$. To do so, we follow the same argument used to obtain \eqref{a10} except for the fact that, defining $\beta_1:=\frac{9^+}8+\alpha_1$, $\partial_x\rho_{N,0}^{\beta_1}$ is no longer bounded above independent of $N$; that is, estimate \eqref{a3} is no longer valid.\\

Instead of \eqref{a3} we argue as follows
\begin{equation}
    \label{b3}
    \begin{split}
        \left| \int_0^{t_1}\int v^4_m \partial_x\rho_{N,0}^{\beta_1}(x)dxdydt\right|&\le c\int_0^{t_1}\|\langle x\rangle^{1/2}v_m\|_{L^\infty_{xy}}^2\int v_m^2 \langle x\rangle dxdydt\\&\le c(T)\|\langle x\rangle^{1/2}v_m\|^2_{L^3_TL^\infty_{xy}}\|\langle x\rangle^{1/2}v_m\|_{L^\infty_TL^2_{xy}}^2\\&\le c(T)\|\langle x\rangle^{\frac{9}{8}^+}v_m\|^2_{L^3_TL^\infty_{xy}}\|\langle x\rangle^{\frac{9}{8}^+}v_m\|_{L^\infty_TL^2_{xy}}^2\\&\le c(T,\eta_{9/4+}(v),\|v\|_{Z_{\frac{9}{4}^+,\frac{9}{8}^+}}),
    \end{split}
\end{equation}
where we have used Lemma \ref{remark2}. 

Taking into account \eqref{b3} and repeating the steps \eqref{a1}-\eqref{a9} done in the proof of Theorem \ref{teo1}, it can be seen that 
\begin{equation}
    \label{b4} \langle x\rangle^{\beta_1-1/2}\partial_xv(t)\in L^2(\R^2)\ a. e. \ t\in[0,t_1].
\end{equation}
Moreover, if the convergence in $m$ is omitted, it follows for $m$ large enough that
\begin{equation}
    \label{b5} \|\langle x\rangle^{\beta_1-1/2}\partial_x v_m\|_{L^2_{xyt_1}}\le c(T,\eta_{9/4+}(v),\|v\|_{Z_{\frac{9}{4}^+,\frac{9}{8}^+}}).
\end{equation}

Without loss of generality we can assume \eqref{b4} holds for $0$ and $t_1$, otherwise we would consider a subinterval $[t_0^*,t_1^*]\subset [0,t_1]$ with $t_0^*$ and $t_1^*$ satisfying \eqref{b4}.\\

Let us take the simmetrized mZK equation \eqref{smZK} for $v_m$ and differentiate once with respect to $x$. Denoting $w_m:=\partial_x v_m$ we get 
\begin{equation}
    \label{b1}\partial_t w_m + \partial_x^3 w_m+\partial_y^3 w_m + \mu \partial_x(v_m^2w_m+v^2_m \partial_y v_m)=0.
\end{equation}
For $i=1,2$ we define $\beta_i:=\frac{9}{8}^++\min\{\alpha_i,\frac{3}{8}^-\}$. Let us consider $\rho_N(x)\equiv \rho_{N,1}^{\beta_1}(x)$. Multiplying \eqref{b1} by $w_m\rho_N(x)$ and integrating by parts we get
\begin{equation}\label{b2}
\begin{split}
    \frac{1}{2}\frac{d}{dt}\int w_m^2\rho_Ndxdy+\frac{3}{2}\int (\partial_x w_m)^2\rho_N'dxdy&-\frac{1}{2}\int w_m^2\rho_N^{\prime \prime \prime}dxdy\\&+\mu\int\partial_x(v^2_m(w_m+\partial_yv_m))w_m\rho_Ndxdy=0,
\end{split}
\end{equation}

where, integrating by parts,
\begin{align}
\notag\int \partial_x(v_m^2(w_m+\partial_yv_m))w_m\rho_N dx dy =& \int v_m^2\partial_x w_m w_m \rho_N dx dy + \int  2v_mw_m^3\rho_N dx dy \\
\notag&+ \int v_m^2\partial_yw_mw_m\rho_N dx dy + \int 2v_mw_m^2\partial_yv_m\rho_N dx dy\\
\notag=& -\int v_m^2\partial_x w_m w_m \rho_N dx dy - \int  v_m^2 w_m^2 \rho'_N dx dy \\
\notag&+ \int v_m^2\partial_yw_mw_m\rho_N dx dy + \int 2v_mw_m^2\partial_yv_m\rho_N dx dy\\
=&I+II+III+IV. \label{nl_int_by_p}
\end{align}

Let us observe that 
\begin{equation}
    \label{b6}
    \begin{split}
        \left|\int_0^{t_1} \int w_m^2 \rho_N^{\prime\prime\prime}dxdydt \right|&\le c\int_0^{t_1}\int w_m^2dxdydt\le c(T)\|w_m\|_{L^\infty_TL^2_{xy}}\le c(T,\eta_1(v)).
    \end{split}
\end{equation}

Now, let us integrate \eqref{nl_int_by_p} over $[0,t_1]$, and estimate each of the four integrals that result on the right hand side.\\

By using Lemma \ref{remark1}, \eqref{b5} and the continuous dependence, we begin estimating
\begin{equation}
    \label{b7}
    \begin{split}
      \left| \int_0^{t_1} I dt\right|  &=\left|\int_0^{t_1}\int v_m^2\partial_x w_m w_m \rho_Ndxdydt \right|\le \|v_m^2w_m\partial_xw_m\rho_N\|_{L^1_{xyt_1}}\\&\le \|\rho_N^{\frac{3}{4}^+}v_m^2\|_{L^2_xL^\infty_{yt_1}}\|\partial_x w_m\|_{L^\infty_xL^2_{yt_1}}\|\rho_N^{\frac{1}{4}^-}w_m\|_{L^2_{xyt_1}}\\&\le c\|\langle x\rangle^{\frac{3}{4}^+}v_m\|_{L^4_xL^\infty_{yt_1}}^2\|\partial_x w_m\|_{L^\infty_xL^2_{yt_1}}\|\langle x\rangle^{1/2}w_m\|_{L^2_{xyt_1}}\\&\le c(T,\eta_1(v), \eta_{9/4+}(v),\|v\|_{Z_{\frac{9}{4}^+,\frac{9}{8}^+}}),
    \end{split}
\end{equation}

We continue with
\begin{equation}
    \label{b9}\begin{split}
        \left|\int_0^{t_1} II dt \right|&=\left|\int_0^{t_1}\int v_m^2w_m^2\rho_N'dxdydt \right|\le c\int_0^{t_1}\|v_m\|^2_{L^\infty_{xy}}\int w_m^2dxdydt\\
        &\le c(T)\|w_m\|_{L^{\infty}_TL^2_{xy}}^2\|v_m\|_{L^3_TL^\infty_{xy}}^2\le c(T,\eta_1(v),\eta_{9/4+}(v)).
    \end{split}
\end{equation}

For the third term we have
\begin{equation}
    \label{b10}
    \begin{split}
       \left|\int_0^{t_1} III dt \right| &= \left|\int_0^{t_1}\int v_m^2\partial_x\partial_y v_m w_m \rho_Ndxdydt \right|\le \|v_m^2\partial_x\partial_yv_mw_m\rho_N\|_{L^1_{xyt_1}}\\&\le\|v_m^2\rho_{N}^{\frac{3}{4}^+}\|_{L^2_xL^\infty_{yt_1}}\|\partial_yw_m\|_{L^\infty_xL^2_{yt_1}}\|\rho_N^{\frac{1}{4}^-}w_m\|_{L^2_{xyt_1}}\le c(T,\eta_1(v),\eta_{9/4+}(v),\|v\|_{Z_{\frac{9}{4}^+,\frac{9}{8}^+}}).
    \end{split}
\end{equation}

For the fourth term we argue as follows
\begin{equation}
    \label{b11}\begin{split}
       \left| \int_0^{t_1} IV dt \right| &= \left|\int_0^{t_1}\int v_m w_m^2\partial_y v_m \rho_N dxdydt \right| \le \|v_m\|_{L^{\infty}_{xyt_1}}\|\partial_y v_m\|_{L^{9/4}_TL^{\infty}_{xy}}\|w_m^2\rho_N\|_{L^{9/5}_{t_1}L^1_{xy}}\\&\le c(T)\|v_m\|_{L^\infty_TH^{9/4^+}_{xy}}\|\partial_yv_m \|_{L^{9/4}_TL^{\infty}_{xy}}\|\rho_N^{1/2}w_m\|^2_{L^2_{xyt_1}}\le c(T,\eta_{9/4+}(v),\|v\|_{Z_{\frac{9}{4}^+,\frac{9}{8}^+}}).
    \end{split}
\end{equation}

From \eqref{nl_int_by_p} to \eqref{b11} and integrating \eqref{b2} over $[0,t_1]$ we conclude
\begin{equation}
    \label{b12}
    \begin{split}
        \int_0^{t_1}&\int (\partial_x w_m)^2\rho_N'dxdydt\\
        &\le c(T,\eta_{1}(v),\eta_{9/4+}(v),\|v\|_{Z_{\frac{9}{4}^+,\frac{9}{8}^+}})+\|\rho_N^{1/2}\partial_xv_m(t_1)\|_{L^2_{xy}}+\|\rho_N^{1/2}\partial_xv_{0m}\|_{L^2_{xy}}\\&\le c(T,\eta_{1}(v),\eta_{9/4+}(v),\|v\|_{Z_{\frac{9}{4}^+,\frac{9}{8}^+}})+c\|\rho^{1/2}_N\partial_xv(t_1)\|_{L^2_{xy}}+c\|\rho^{1/2}_N\partial_xv_0\|_{L^2_{xy}}\\&\le c(T,\eta_{1}(v),\eta_{9/4+}(v),\|v\|_{Z_{\frac{9}{4}^+,\frac{9}{8}^+}})+c\|\langle x\rangle^{\beta_1-1/2}\partial_xv(t_1)\|_{L^2_{xy}}+c\|\langle x\rangle^{\beta_1-1/2}\partial_xv_0\|_{L^2_{xy}}\\&\le c(T,\eta_{1}(v),\eta_{9/4+}(v),\|v\|_{Z_{\frac{9}{4}^+,\frac{9}{8}^+}}, \|\langle x\rangle^{\beta_1-1/2}\partial_xv(t_1)\|_{L^2_{xy}},\|\langle x\rangle^{\beta_1-1/2}\partial_xv_0\|_{L^2_{xy}}).
    \end{split}
\end{equation}

From \eqref{b12}, following the same ideas of the proof of Theorem \ref{teo1} it can be seen that 
\begin{equation}\label{b13}
    \begin{cases}
        \langle x \rangle^{\beta_1-1}\partial_x^2 v(t)\in L^2(\R^2) \ a.e.\ t\in[0,t_1],\\
        \langle y\rangle^{\beta_2-1}\partial_y^2 v(t)\in L^2(\R^2) \ a.e.\ t\in[0,t_1]\ \mbox{(by symmetry).}
    \end{cases}
\end{equation}

Without loss of generality assume \eqref{b13} holds for $0$ and $t_1$.\\

Now, we differentiate the simmetrized mZK equation \eqref{smZK} for $v_m$ twice respect to the variable $x$, multiply by $\partial_x^2 v_m \rho_N(x)$, where $\rho_N(x):=\rho_{N,2}^{\beta_1}$, use integration by parts, and denote $z_m:=\partial_x^2 v_m$, to obtain
\begin{equation}\label{b14}
    \begin{split}
        \frac{1}{2}\frac{d}{dt}\int z_m^2\rho_N dxdy +\frac{3}{2}\int (\partial_xz_m)^2&\rho_N'dxdy-\frac{1}{2}\int z_m^2\rho_N^{\prime\prime\prime}dxdy\\&+\mu\int\partial_x^2(v_m^2(\partial_xv_m+\partial_yv_m))z_m\rho_Ndxdy=0.
    \end{split}
\end{equation}

Let us note that the linear part of \eqref{b14} is analogous to the linear part of \eqref{b2} when seen for $z_m$. We therefore focus on the nonlinear part. 

Note \begin{equation}
\begin{split}
    z_m\rho_N\partial_x^2(v_m^2(\partial_x v_m+\partial_yv_m))&=z_m\rho_N\sum_{j=0}^{2}{2\choose j}\partial_x^j(v_m^2)(\partial_x^{3-j}v_m  + \partial_x^{2-j}\partial_yv_m)\\&\hspace{-30mm}=\sum_{j=0}^2\sum_{k=0}^j{2 \choose j}{j \choose k}z_m\rho_N\left\{\partial_x^{3-j}v_m\partial_x^{j-k}v_m \partial_x^k v_m + \partial_x^{2-j}\partial_y v_m \partial_x^{j-k}v_m \partial_x^{k}v_m \right\}\\&\hspace{-30mm}\equiv\sum_{j=0}^2\sum_{k=0}^j z_m\rho_N\left\{A_{jk}+B_{jk} \right\}.
\end{split}
\end{equation}

Let us first estimate $$I_{jk}\equiv \left|\int_0^{t_1}\int z_m\rho_N A_{jk}dxdydt \right|.$$ 

\noindent\underline{If j=0}:\\
\begin{equation}
    \label{b15}
    \begin{split}
        I_{00}&\le \|\partial_x^3v_m\|_{L^\infty_xL^2_{yt_1}}\|\rho_N^{1/2}v_m\|_{L^4_xL^{\infty}_{yt_1}}^2\|\partial_x^2v_m\|_{L^2_{xyt_1}}\\&\le c(T)\|\partial_x^2\partial_xv_m\|_{L^\infty_xL^2_{yt_1}}\|\langle x\rangle^{\frac{3}{4}^+}v_m\|_{L^4_xL^\infty_{yt_1}}^2\|z_m\|_{L^\infty_TL^2_{xy}}\le c(T,\eta_{9/4+}(v),\|v\|_{Z_{\frac{9}{4}^+,\frac{9}{8}^+}}).
    \end{split}
\end{equation}

\noindent\underline{If j=2}: First, in case $k\neq 1$ we have, using Lemma \ref{remark2}, that \begin{equation}
    \label{b16}
    \begin{split}
        I_{22}=I_{20}&\le \int_0^{t_1}\|\rho_Nv_m\|_{L^\infty_{xy}}\|\partial_x v_m\|_{L^\infty_{xy}}\int z_m^2dxdydt \\&\le c(T)\|\partial_x v_m\|_{L^\infty_{xyt_1}}\|\rho_N v_m\|_{L^3_TL^\infty_{xy}}\|z_m\|_{L^\infty_TL^2_{xy}}^2\le c(T,\eta_2(v),\|v\|_{Z_{\frac{9}{4}^+,\frac{9}{8}^+}}).
    \end{split}
\end{equation}
In case $k=1$ we argue as follows
\begin{equation}
    \label{b17}
    \begin{split}
        I_{2,1}&\le \|\partial_x v_m\|_{L^{\infty}_{xyt_1}}^2\|\rho_N\partial_x v_m\|_{L^2_{xyt_1}}\|\partial_x^2 v_m\|_{L^2_{xyt_1}}\\&\le c(T)\|v_m\|^3_{L^\infty_TH^{9/4^+}_{xy}}\|\langle x\rangle^{2\beta_1-2}\partial_x v_m\|_{L^2_{xyt_1}}\le c(T,\vertiii{v})\|\langle x\rangle^{\beta_1-1/2}\partial_x v_m\|_{L^2_{xyt_1}}\\&\le c(T,\eta_{9/4+}(v),\|v\|_{Z_{\frac{9}{4}^+,\frac{9}{8}^+}}),
    \end{split}
\end{equation}
where we have used \eqref{b5}.\\

\noindent\underline{If j=1}: \begin{equation}
    \label{b18}I_{1k}=I_{22}\le c(T,\eta_2(v),\|v\|_{Z_{\frac{9}{4}^+,\frac{9}{8}^+}}), \text{ for }k=0,1.
\end{equation}

On the other hand in what comes to the $B_{jk}$ terms, we proceed to estimate $$\RNum{2}_{jk}\equiv \left|\int_0^{t_1}\int z_m\rho_N B_{jk}dxdydt \right|$$ as follows.\\

\noindent\underline{If j=0}: Same estimate \eqref{b15} holds with $\|\partial_{x}^2\partial_y v_m\|_{L^\infty_xL^2_{yt_1}}$ instead of $\|\partial_x^3v_m\|_{L^\infty_{x}L^2_{yt_1}}$. \\

\noindent\underline{If j=2}: When $k\neq 1$ the estimate is analogous to \eqref{b16} and to \eqref{b17} when $k=1$.\\

\noindent\underline{If j=1}: We have that $II_{10}=II_{11}$. Besides, by using integration by parts, we have that
\begin{equation}\label{int_parts_1}
    \begin{split}
        \int\partial_x\partial_yv_m v_m \partial_xv_mz_m\rho_Ndxdy&=-\int \partial_y v_m (\partial_xv_m)^2z_m\rho_Ndxdy-\int \partial_yv_m v_m z_m^2\rho_Ndxdy\\&\hspace{-20mm} -\int \partial_yv_m v_m \partial_x v_m \partial_x^3v_m\rho_Ndxdy -\int \partial_y v_m v_m \partial_x v_m \partial_x^2 v_m \rho_N'dxdy.
    \end{split}
\end{equation}

Taking into account that estimate \eqref{b17} holds with $\|\partial_x v_m\|_{L^\infty_{xyt_1}} \|\partial_y v_m\|_{L^\infty_{xyt_1}}$ instead of $\|\partial_x v_m\|^2_{L^\infty_{xyt_1}}$, estimate \eqref{b16} holds with $\|\partial_y v_m\|_{L^\infty_{xyt_1}}$ instead of $\|\partial_x v_m\|_{L^\infty_{xyt_1}}$, and the facts
\begin{align*}
\left|\int_0^{t_1} \int \partial_y v_m v_m \partial_x v_m \partial_x^2 v_m \rho_N dx dy dt \right| & \leq \|\partial_x^2 \partial_x v_m\|_{L^\infty_x L^2_{yt_1}} \|\partial_y v_m\|_{L^\infty_{xyt_1}}\|\rho_N \partial_x v_m\|_{L^2_{xyt_1}} \|v_m\|_{L^2_x L^\infty_{yt_1}}\\
&\leq c(\eta_2(v),\eta_{9/4+}(v)),
\end{align*}

and
\begin{equation*}\begin{split}
    \left|\int_0^{t_1} \int \partial_y v_m v_m \partial_xv_m \partial_x^2 v_m \rho_N'dxdydt \right|&\le c(T)\|\partial_y v_m\|_{L^\infty_{xyt_1}}\|\partial_x v_m\|_{L^\infty_{xyt_1}}\|v_m\|_{L^\infty_TL^2_{xy}}\|\partial_x^2 v_m\|_{L^\infty_TL^2_{xy}}\\&\le c(T,\eta_2(v)),
    \end{split}
\end{equation*}
 we conclude, integrating \eqref{int_parts_1} over $[0,t_1]$, that
\begin{equation}
    \label{b19} \RNum{2}_{1k}\le c(T,\eta_2(v),\eta_{9/4+}(v),\|v\|_{Z_{\frac{9}{4}^+,\frac{9}{8}^+}}) \text{ for }k=1,2,
\end{equation}

Integrating \eqref{b14} over $[0,t_1]$ and taking into account \eqref{b13}-\eqref{b19} it can be seen as done in \eqref{b12} that
\begin{equation}
    \label{b20} \int_0^{t_1}\int (\partial_x z_m)^2\rho_Ndxdydt \le c(T,\eta_2(v),\eta_{9/4+}(v),\|v\|_{Z_{\frac{9}{4}^+,\frac{9}{8}^+}}, \|\langle x\rangle^{\beta_1-1}\partial_xv(t_1)\|_{L^2_{xy}},\|\langle x\rangle^{\beta_1-1}\partial_xv_0\|_{L^2_{xy}}).
\end{equation}

From \eqref{b20} following the ideas of the proof of Theorem \ref{teo1} it can be seen that 
\begin{equation}\label{b21}
    \begin{cases}
        \langle x\rangle^{\beta_1-3/2}\partial_x^3v(t) \in L^2(\R^2)\ a.e.\ t\in[0,t_1],\\
        \langle y\rangle^{\beta_2-3/2}\partial_y^3v(t) \in L^2(\R^2)\ a.e.\ t\in[0,t_1]\ \mbox{(by symmetry),}
    \end{cases}
\end{equation}
where $\beta_2:=\frac{9^+}8+\alpha_2$.\\

Let $t^*\in[0,t_1]$ be such that \eqref{b13} and \eqref{b21} hold for $t^*$. For $i=1,2$ we define $f_i(t^*)(x,y):=\langle a_i \rangle^{\beta_i-3/2}\partial_x^2v(t^*)(x,y)$ where $a_1=x$ and $a_2=y$. Interpolating with $\theta=2\beta_i-2$ we get
\begin{equation}
    \label{b22} \|J^{2\beta_i-2}_{a_i}\partial_{a_i}^2v(t^*)\|_{L^2_{xy}}=\|J^{\theta\cdot 1}_{a_i} \langle a_i\rangle^{(1-\theta)\cdot\frac{1}{2}}f_i \|_{L^2_{xy}}\le c \|\langle a_i\rangle^{1/2}f_i(t^*)\|_{L^2_{xy}}^{1-\theta}\|J^{1}_{a_i}f_i(t^*)\|_{L^2_{xy}}^\theta<\infty.
\end{equation}

From \eqref{b22} we conclude $v(t^*)\in H^{\frac{9}{4}^++2\alpha}(\R^2)$ which implies $$v\in C([0,T];H^{\frac{9}{4}^++2\alpha}(\R^2)).$$ Moreover, from Theorem \ref{LWPZ}
\begin{equation}\label{b23}
    v\in C([0,T];Z_{\frac{9}{4}^++2\alpha,\frac{9}{8}^++\alpha}).
\end{equation}

\noindent \underline{\textbf{Case $\alpha\in [\frac{4\eta-1}{8},\frac{4\eta+3}{8})$, $\eta\in\N$}.}\\

Let us consider the following statement \\

\noindent\textbf{P(n)}: For $i=1,2$ and all $\beta_i\in\left(\frac{4n+8}{8},\frac{4n+12}{8}\right]$ with $ \beta_i\le\frac{9}{8}^++\min\left\{\alpha_i,\frac{4n+3}{8}^-\right\}$ there is a non-degenerated subinterval $I\subset[0,t_1]$ such that for almost every $t\in I$ we have 
\begin{equation}\label{hyp}
    \begin{cases}
        \langle a_i\rangle^{\beta_i-1-n/2}\partial_{a_i}^{n+2}v(t) \in L^2(\R^2),\\
        \langle a_i \rangle^{\beta_i-3/2-n/2}\partial_{a_i}^{n+3}v(t)\in L^2(\R^2);
    \end{cases}
\end{equation}
where $a_1=x$ and $a_2=y$. \\

It is enough to use induction to prove P(n) for all $n\in\{0,1,\dots,\eta\}$. In such case, taking $t^*\in[0,t_1]$ such that \eqref{hyp} holds for $n=\eta$ with $\beta_i=\frac{9}{8}^++\min\left\{\alpha_i,\frac{4\eta+3}{8}^-\right\}$ and considering $f_i(x,y):=\langle a_i\rangle^{\beta_i-3/2-\eta/2}\partial_{a_i}^{\eta+2}v(t^*)(x,y)$ we would interpolate with $\theta=2\beta_i-2-\eta$ to get \begin{equation}
    \|J^{2\beta_i-2-\eta}_{a_i}\partial_{a_i}^{\eta+2}v(t^*) \|_{L^2_{xy}}<\infty,
\end{equation}
which implies $v\in C([0,T];H^{\frac{9}{4}^++2\alpha}(\R^2))$.\\

The base case ${P(0)}$ follows from \eqref{b13} and \eqref{b21}. \\

Let us assume that ${P(n)}$ is true and let us prove ${P(n+1)}$. Without loss of generality we assume $I=[0,t_1]$. From \eqref{hyp} with $\beta_1=\beta_2=\frac{4n+12}{8}=\frac{n+3}{2}$ it follows that $\partial_{a_i}^{n+3}v(t)\in L^2(\R^2)$ $a.e.\ t\in[0,t_1]$; which implies $v\in C([0,T];H^{n+3}(\R^2))$. Since $n<\eta$ we have $\frac{3+n}{2}\le \frac{9}{8}^++\alpha_i$. Therefore, using Theorem \ref{LWPZ} we conclude $v\in C([0,T];Z_{n+3,\frac{n+3}{2}})$. 

Note from Lemma \ref{hilberpp} and \eqref{interpol1r} that
\begin{equation}
    \label{A3} \|\langle x\rangle^{\frac{n+3-k}{2}}\partial_x^k v \|_{L^2_{xy}}\le c\|J^{n+3}_{xy} v\|_{L^2_{xy}}^\theta \|\langle x\rangle^{\frac{n+3}{2}}v\|_{L^2_{xy}}^{1-\theta},\ \ \theta=\frac{n+3-k}{n+3}.
\end{equation}
Also, from Sobolev embeddings, for $k\in\{0,1,\dots,n+1 \}$
\begin{equation}
    \label{A4} \|\partial_x^k v\|_{L^\infty_{xy}}\le c\|v\|_{H^{n+3}_{xy}}.
\end{equation}

Let $\{v_{0m}\}\subset C_0^{\infty}(\R^2)$ be a sequence converging to $v_0$ in $Z_{n+3,\frac{n+3}{2}}$. Denote with $v_m$ the solution provided by Theorem \ref{LWPZ} with initial datum $v_{0m}$. Assume for all $m$ large enough that $v_m\in C([0,T];Z_{n+3,\frac{n+3}{2}})$. 

Let $\beta_1\in \left(\frac{4n+12}{8},\frac{4n+16}{8} \right]=\left(\frac{n+3}{2},\frac{n+4}{2} \right]$ with $\beta_1\le \frac{9}{8}^++\min\{\alpha_i,\frac{4\eta+3}{8}^-\}$.

Let us denote $w_m:=\partial_x^{n+2}v_m$. Differentiating the symmetrized mZK equation \eqref{smZK} for $v_m$, $n+2$ times respect to the variable $x$, we get
\begin{equation}
    \label{A5}\partial_t w_m + \partial_x^3 w_m +\partial_y^3w_m+\mu\partial_x^{n+2}(v_m^2\partial_xv_m+v_m^2\partial_yv_m)=0.
\end{equation}
Let us multiply \eqref{A5} by $w_m \rho_n(x)$, where $\rho_N(x):= \rho_{N,n+2}^{\beta_1}$, and integrate by parts to get
\begin{equation}
    \label{A6}\begin{split}
        \frac{1}{2}\frac{d}{dt}\int w_m^2\rho_N dxdy+\frac{3}{2}\int (\partial_x w_m)^2&\rho_N'dxdy-\frac{1}{2}\int w_m^2\rho_N^{\prime\prime\prime}dxdy\\&+\mu\int\partial_x^{n+2}(v_m^2\partial_xv_m+v_m^2\partial_y v_m)w_m\rho_Ndxdy=0.
    \end{split}
\end{equation}

We note that the linear part of \eqref{A5} is the same as in the symmetrized mZK equation when seen for $w_m$, therefore to argue as in the proof of Theorem \ref{teo1} we only need to estimate $$\left|\int_0^{t_1}\int \partial_x^{n+2}(v_m^2 \partial_x v_m +v_m^2\partial_yv_m)w_m\rho_Ndxdtdy \right|.$$
Since $$\partial_x^{n+2}(v_m^2 \partial_x v_m +v_m^2\partial_yv_m)=\sum_{j=0}^{n+2}\sum_{k=0}^j{n+2 \choose j}{j \choose k}\partial_x^{j-k}v_m\partial_x^kv_m(\partial_x^{n+3-j}v_m+\partial_x^{n+2-j}\partial_y v_m),$$
we have $$\left|\int_0^{t_1}\int \partial_x^{n+2}(v_m^2 \partial_x v_m +v_m^2\partial_yv_m)w_m\rho_Ndxdtdy \right|\le \sum_{j=0}^{n+2}\sum_{k=0}^j{n+2\choose j}{j\choose k}(A_{jk}+B_{jk}),$$ where $$
    A_{jk}:=\|\partial_x^{n+3-j}v_m\partial_x^{j-k}v_m \partial_x^k v_m \partial_x^{n+2}v_m\rho_N\|_{L^1_{xyt_1}}$$
and $$B_{jk}:=\|\partial_x^{n+2-j}\partial_yv_m\partial_x^{j-k}v_m \partial_x^k v_m \partial_x^{n+2}v_m\rho_N\|_{L^1_{xyt_1}}.$$

We first estimate the $A_{jk}$ terms.\\

\noindent\underline{If $j=0$}:\\
\begin{equation}
    \label{A7} \begin{split}
        \|\partial_x^{n+3}v_m v_m^2 \partial_x^{n+2}v_m\rho_N\|_{L^1_{xyt_1}}&\le \|\rho_N v_m^2\|_{L^1_TL^\infty_{xy}}\|\partial_x^{n+3}v_m\|_{L^\infty_TL^2_{xy}}\|\partial_x^{n+2}v_m\|_{L^\infty_TL^2_{xy}}\\&\le c(T)\|\rho_N^{1/2}v_m\|^2_{L^3_TL^\infty_{xy}}\|v_m\|^2_{L^\infty_TH^{n+3}_{xy}} \\& \le c(T)\|\langle x \rangle v_m\|^{2}_{L^3_TL^\infty_{xy}}\|v_m\|_{L^\infty_T H^{n+3}_{xy}} \\&\le c(T,\eta_2(v),\|v\|_{Z_{n+3,\frac{n+3}{2}}}), 
    \end{split}
\end{equation}
where we have used Lemma \ref{remark2}. \\

\noindent\underline{If $j=n+2$}:\\

First, if $k=0$ or $k=n+2$ we have, by using Lemma \ref{remark2}, that
\begin{equation}
    \label{A8}\begin{split}
        \|\partial_xv_m (\partial_x^{n+2}v_m)^2v_m\rho_N\|_{L^1_{xyt_1}}&\le \int_0^{t_1}\|\partial_x v_m\|_{L^\infty_{xy}} \|\rho_Nv_m\|_{L^\infty_{xy}}\int (\partial_x^{n+2}v_m)^2dxdy\\&\le c(T)\|v_m\|^3_{L^\infty_TH^{n+3}_{xy}}\|\langle x\rangle^2 v_m\|_{L^3_TL^\infty_{xy}}\le c(T,\|v\|_{Z_{n+3,\frac{n+3}{2}}}).
    \end{split}
\end{equation}

Now, if $k\in\{1,2,\dots,n+1\}$, we have that
\begin{equation}\label{A9}
\begin{split}
    \|\partial_x v_m \partial_x^{n+2-k}&v_m\partial_x^kv_m \partial_x^{n+2}v_m \rho_N\|_{L^1_{xyt_1}}\\&\le c(T)\|\rho_N^{3/4}\partial_xv_m\|_{L^\infty_TL^2_{xy}}\|\rho^{1/4}_N\partial_x^{n+2}v_m\|_{L^2_{xyt_1}}\|\partial_x^k v_m\|_{L^\infty_{xyt_1}}\|\partial_x^{n+2-k}v_m\|_{L^\infty_{xyt_1}}\\&\le c(T)\|\langle x\rangle^{3/2}\partial_x v_m\|_{L^\infty_TL^2_{xy}}\|\langle x\rangle^{1/2}\partial_x^{n+2}v_m\|_{L^2_{xyt_1}}\|v_m\|_{L^\infty_TH^{n+3}_{xy}}^2\\&\le c(T)\|\langle x \rangle^{\frac{n+3-1}{2}}\partial_x v_m\|_{L^\infty_TL^2_{xy}}\|\langle x\rangle^{1/2}\partial_x^{n+2}v_m\|_{L^2_{xyt_1}}\|v_m\|_{L^\infty_TH^{n+3}_{xy}}^2\\&\le c(T,\|v\|_{Z_{n+3,\frac{n+3}{2}}}),
\end{split}
\end{equation}
where we have used \eqref{A3}, the inductive hypothesis, and \eqref{A4}.\\ 

\noindent\underline{If $j=1$}: In this case we follow the exact same estimate \eqref{A8}.\\

\noindent\underline{If $j=2$}:\\

First, when $k=0$ or $k=2$ we have \begin{equation}
    \label{A10}\begin{split}
        \|\partial_x^{n+1} v_m &\partial^2_xv_m v_m \partial_x^{n+2}v_m \rho_N \|_{L^1_{xyt_1}}\\&\le c(T)\|\partial_x^{n+1}v_m\|_{L^\infty_{xyt_1}} \|\partial^2_xv_m\|_{L^\infty_{xyt_1}} \|\langle x\rangle^{2}v_m\|_{L^\infty_TL^2_{xy}}\|\partial_x^{n+2}v_m\|_{L^\infty_TL^2_{xy}}\\&\le c(T)\|v_m\|^3_{L^\infty_TH^{n+3}_{xy}}\|\langle x \rangle^{\frac{n+3}{2}}v_m\|_{L^\infty_TL^2_{xy}}\\&\le c(T,\|v\|_{Z_{n+3,\frac{n+3}{2}}}).
    \end{split}
\end{equation}
Now, in case $k=1$ we have
\begin{equation}
    \label{A11}\begin{split}
        \|\partial_x^{n+1}v_m &(\partial_x v_m)^2\partial_x^{n+2}v_m\rho_N\|_{L^1_{xyt_1}}\\& \le c(T)\|\partial_x^{n+1}v_m\|_{L^\infty_{xyt_1}}\|\rho^{3/4}_N\partial_xv_m\|_{L^\infty_TL^2_{xy}}\|\rho_N^{1/4}\partial_x^{n+2}v_m\|_{L^2_{xyt_1}}\|\partial_xv_m\|_{L^\infty_{xyt_1}}\\&\le c(T) \|v_m\|^2_{L^\infty_TH^{n+3}_{xy}}\|\langle x\rangle^{\frac{n+3-1}{2}}\partial_xv_m\|_{L^\infty_TL^2_{xy}}\|\langle x\rangle^{1/2}\partial_x^{n+2}v_m\|_{L^2_{xyt_1}}\\&\le c(T,\|v\|_{Z_{n+3,\frac{n+3}{2}}}),
    \end{split}
\end{equation}
where we have used the inductive hypothesis.\\

\noindent\underline{If $j\in\{3,4,\dots,n+1\}$}:\\

\begin{equation}
    \label{A12}\begin{split}
        \|\partial_x^{n+3-j}v_m\partial_x^{j-k}&v_m \partial_x^kv_m\partial_x^{n+2}v_m \rho_N \|_{L^1_{xyt_1}}\\&\le c(T)\|\partial_x^{j-k}v_m\|_{L^\infty_{xyt_1}}\|\partial_x^kv_m\|_{L^\infty_{xyt_1}}\|\rho_N^{3/4}\partial_x^{n+3-j}v_m\|_{L^\infty_TL^2_{xy}}\|\rho^{1/4}_N \partial_x^{n+2}v_m\|_{L^2_{xyt_1}}\\&\le c(T)\|v_m\|^2_{L^\infty_TH^{n+3}_{xy}} \|\langle x \rangle^{j/2}\partial_x^{n+3-j}v_m\|_{L^\infty_TL^2_{xy}}\|\langle x\rangle^{1/2}\partial_x^{n+2}v_m\|_{L^2_{xyt_1}}\\&\le c(T,\|v\|_{Z_{n+3,\frac{n+3}{2}}}),
    \end{split}
\end{equation}
where we used \eqref{A3} and the inductive hypothesis.\\

We now proceed to estimate the $B_{kj}$ terms. \\

\noindent\underline{If $j=0$}: Analogous to \eqref{A7} with $\partial_x^{n+2} \partial_y v_m$, instead of $\partial_x^{n+3}v_m$.\\

\noindent\underline{If $j=n+2$}: When $k=0$ or $k=n+2$ we argue analogous to \eqref{A8} with $\partial_y v_m$, instead of $\partial_x v_m$. When $k\in\{1,2,\dots,n+1\}$ the estimate follows analogous to \eqref{A9} with $\partial_y v_m$ instead of $\partial_x v_m$, using the following interpolation \begin{equation}
    \label{A13}
    \begin{split}
        \|\rho^{3/4}_N \partial_y v_m\|_{L^\infty_TL^2_{xy}}&\le \|\langle x\rangle^{3/2}\partial_yv_m\|_{L^\infty_TL^2_{xy}}=\|\partial_y(\langle x\rangle^{3/2}v_m)\|_{L^\infty_TL^2_{xy}}\le \|J^1_{xy}\langle x\rangle^{3/2}v_m\|_{L^\infty_TL^2_{xy}}\\&\le c \|J^{n+3}_{xy}v_m\|_{L^\infty_TL^2_{xy}}^\theta \|\langle x \rangle^{\frac{n+3}{2}}v_m\|_{L^\infty_TL^2_{xy}}^{1-\theta},
    \end{split}
\end{equation}
where $\theta=\frac{1}{n+3}$, $a=n+3$ and $b=\frac{3(n+3)}{2(n+2)}\le \frac{n+3}{2}$.\\

\noindent\underline{If $j=1$}: Analogous to \eqref{A8} with $\partial_x^{n+1}\partial_y v_m \partial_x^{n+2} v_m$ instead of $(\partial_x^{n+2} v_m)^2$.\\

\noindent\underline{If $j=2$}: Analogous to \eqref{A10} and \eqref{A11} with $\partial_x^n \partial_y v_m$ instead of $\partial_x^{n+1} v_m$.\\

\noindent\underline{If $j\in\{3,4,\cdots,n+1\}$}: Analogous to \eqref{A12} with $\partial_x^{n+2-j}\partial_y v_m$ instead of $\partial_x^{n+3-j}v_m$, taking into account that for $\theta=\frac{j}{n+2}$
\begin{equation}
    \|\rho_N^{3/4}\partial_x^{n+2-j}\partial_yv_m\|_{L^2_{xy}}\le c\|\langle x\rangle^{3/2}\partial_x^{n+2-j}\partial_y v_m\|_{L^2_{xy}}\le c\|J^{n+2}_x\partial_yv_m\|_{L^2_{xy}}^{1-\theta}\|\langle x\rangle^{\frac{3n+6}{2j}}\partial_y v_m\|_{L^2_{xy}}^\theta
\end{equation}
and that for $\lambda=\frac{1}{n+3}$ \begin{equation}\label{A15}\begin{split}
    \|\langle x\rangle^{\frac{3n+6}{2j}}\partial_y v_m\|_{L^2_{xy}}&\le c\|J^1_{xy}\langle x\rangle^{\frac{3n+6}{2j}}v_m\|_{L^2_{xy}}\le c\|J^{n+3}_{xy}v_m\|_{L^2_{xy}}^\lambda \|\langle x\rangle^{\frac{3n+9}{2j}}v_m\|_{L^2_{xy}}^{1-\lambda}\\&\le c\|v_m\|_{H^{n+3}_{xy}}+ c\|\langle x\rangle^{\frac{n+3}{2}}v_m\|_{L^2_{xy}}.
    \end{split}
\end{equation}

Taking into account \eqref{A7} - \eqref{A15}, integrating \eqref{A6} over $[0,t_1]$ and arguing as in the proof of Theorem \ref{teo1} it can be seen that
\begin{equation}
    \label{A16}\int_0^{t_1}\int ( \partial_x^{n+3}v)^2\langle x\rangle^{\beta_1-3/2-n/2}dxdydt<\infty,
\end{equation}
which implies \begin{equation}
    \label{A17}\langle x\rangle^{\beta_1-3/2-n/2}\partial_x^{n+3}v(t) \in L^2(\R^2)\ a.e.\ t\in[0,t_1].
\end{equation}

On the other hand, denoting $z_n:= \partial_x^{n+3} v_m$, and differentiating the symmetrized mZK equation \eqref{smZK} $(n+3)$-times with respect to the variable $x$, we get 
\begin{equation}
    \label{A18}
    \partial_tz_m+\partial_x^3z_m +\partial_y^3 z_m +\mu \partial_x^{n+3}(v_m^2\partial_x v_m+v_m^2\partial_y v_m)=0.
\end{equation}

Let us multiply \eqref{A18} by $z_m\rho_N$, where $\rho_N(x):= \rho_{N,n+3}^{\beta_1}(x)$, and integrate by parts to get 
\begin{equation}\label{A19}
    \begin{split}
        \frac{1}{2}\frac{d}{dt}\int z_m^2 \rho_N dxdy 
    +&\frac{3}{2}\int (\partial_x z_m)^2 \rho_N'dxdy
    -  \frac{1}{2}\int z_m^2 \rho_N'''dxdy 
    \\&
     +\mu \int \partial_x^{n+3}(v_m^2\partial_xv_m+v_m^2\partial_yv_m)z_m\rho_Ndxdy=0.
    \end{split}
\end{equation}

To argue as in the proof of Theorem \ref{teo1}, we only need to take care of the nonlinear part. Note that
\begin{align*}
    \left| \int_0^{t_1}\int \partial_x^{n+3}(v_m^2\partial_x v_m+v_m^2\partial_yv_m)z_m\rho_N(x)dxdydt  \right|
     \leq \sum_{j=0}^{n+3}\sum_{k=0}^j \binom{n+3}{j}\binom{j}{k} \{ \bar{A}_{jk}+\bar{B}_{jk} \},
\end{align*}
where 
\begin{equation*}
    \bar{A}_{jk} := \| \partial_x^{n+4-j}v_m\partial_x^{j-k}v_m\partial_x^k v_m z_m\rho_N \|_{L_{xyt_1}^1}
\end{equation*}
and
\begin{equation*}
    \bar{B}_{jk} := \| \partial_x^{n+3-j}\partial_yv_m\partial_x^{j-k}v_m\partial_x^k v_m z_m\rho_N \|_{L_{xyt_1}^1}.
\end{equation*}

We proceed to estimate the $\bar{A}_{jk}$ terms.\\

\noindent\underline{If $j=0$}:\\

\begin{equation}
    \label{A20}
        \begin{split}
            \| \partial_x^{n+4-j}v_m\partial_x^{j-k}v_m\partial_x^k v_m z_m\rho_N \|_{L_{xyt_1}^1} &
            \leq \| \partial_x^{n+4}v_m \|_{L_x^{\infty} L_{yT}^2} \| \rho_N^{\frac{3}{4}^+}v_m^2 \|_{L_x^2L_{yT}^\infty} \| \partial_x^{n+3}v_m \rho_N^{\frac{1}{4}^-} \|_{L_{xyt_1}^2}
            \\&\hspace{-10mm}
             \leq c(T)\| \partial_x^{n+4}v_m \|_{L_x^\infty L_{yT}^2} \| \langle x \rangle^{\frac{3}{4}^+}v_m \|_{L_x^4L_{yT}^\infty}  \| \langle x \rangle^{\beta_1 -3/2-n/2} \partial_x^{n+3}v_m \|_{L_{xyt_1}^2}
            \\&\hspace{-10mm}
            \leq c(T,\eta_{9/4+}(v),\|v\|_{Z_{n+3,\frac{n+3}{2}}}),
        \end{split}
    \end{equation}
where we have used the continuous dependence, Lemma \ref{remark1} and \eqref{A16}.\\

\noindent\underline{If $j=1$}:\\
\begin{equation}
\label{A21}
    \begin{split}
        \| (\partial_x^{n+3}v_m)^2v_m\partial_x v_m\rho_N \|_{L_{xyt_1}^1}&\leq \| v_m \|_{L_{xyt_1}^\infty} \| \partial_xv_m \|_{L_{xyt_1}^\infty} \| \rho_N(\partial_x^{n+3}v_m)^2 \|_{L_{xyt_1}^1}
        \\&
        \leq c \| v_m \|_{L_T^\infty H_{xy}^{n+3}}^2 \|\langle x\rangle^{\beta_1-3/2-n/2}\partial_x^{n+3}v_m\|^2_{L^2_{xyt_1}}
        \\&\leq c(T,\eta_{n+3}(v),\|v\|_{Z_{n+3,\frac{n+3}{2}}}),
    \end{split}
\end{equation}
where we have used \eqref{A16}.\\

\noindent\underline{If $j\in\{2,3,\dots,n+1\}$}:\\

\begin{equation}
\label{A22}
    \begin{split}
        \| \partial_x&^{n+4-j}v_m\partial_x^{j-k}v_m\partial_x^k v_m \partial_x^{n+3}v_m \rho_N \|_{L_{xyt_1}^1}\\&
        \leq c(T)\| \langle x \rangle^{1/2}\partial_x^{n+4-j}v_m \|_{L_T^\infty L_{xy}^2}
        \| \partial_x^{j-k}v_m \|_{L_{xyt_1}^\infty} \| \partial_x^k v_m \|_{L_{xyt_1}^\infty}
        \| \langle x \rangle^{\beta_1 -3/2-n/2}\partial_x^{n+3}v_m \|_{L_{xyt_1}^2}
        \\&        \leq c_T\| v_m \|_{L_T^\infty H_{xy}^{n+3}}^2
        \| \langle x \rangle^{\frac{j-1}{2}}\partial_x^{n+4-j}v_m \|_{L_T^\infty L_{xy}^2}
        \| \langle x \rangle^{\beta_1 -3/2-n/2}\partial_x^{n+3}v_m \|_{L_{xyt_1}^2}
        \\
        &\leq c(T,\eta_{n+3}(v),\|v\|_{Z_{n+3,\frac{n+3}{2}}}),
    \end{split}
\end{equation}
where we have used \eqref{A3} and \eqref{A16}.\\

\noindent\underline{If $j=n+2$}:\\

First, if $k=0$ or $k=n+2$ we have \begin{equation}
        \label{A23}
        \begin{split}            
            \| \partial_x^2 v_m&\partial_x^{n+2-k}v_m \partial_x^kv_m \partial_x^{n+3}v_m \rho_N \|_{L_{xyt_1}^1} 
            \\&
            \leq \| \partial_x^2 v_m \|_{L_{xyt_1}^\infty} \| v_m \|_{L_{xyt_1}^\infty } \| \langle x \rangle^{1/2}\partial_x^{n+2}v_m \|_{L_{xyt_1}^2} \| \langle x \rangle^{\beta_1-3/2-n/2}\partial_x^{n+3}v_m \|_{L_{xyt_1}^2} 
            \\&
            \leq c \| v_m \|_{L_T^\infty H_{xy}^{n+3}}^2 \| \langle x \rangle^{1/2}\partial_x^{n+2}v_m \|_{L_{xyt_1}^2} \| \langle x \rangle^{\beta_1-3/2-n/2}\partial_x^{n+3}v_m \|_{L_{xyt_1}^2} \\&
            \leq c(T,\eta_{n+3}(v),\|v\|_{Z_{n+3,\frac{n+3}{2}}}),
        \end{split}
        \end{equation}
        where we have used \eqref{A16} and the inductive hypotheses. 

Now, If $k\in \{ 1,\dots, n+1 \}$ we have
\begin{equation}
        \label{A24}
            \begin{split}
                \| \partial_x^2 v_m &\partial_x^{n+2-k}v_m \partial_x^k v_m \partial_x^{n+3}v_m \rho_N \|_{L_{xyt_1}^1}
                \\&
                \leq c(T) \| \partial_x^2 v_m \|_{L_{xyt_1}^\infty} \| \partial_x^{n+2-k}v_m \|_{L_{xyt_1}^\infty} \| \langle x \rangle^{1/2}\partial_x^k v_m \|_{L_T^\infty L_{xy}^2} \| \langle x \rangle^{\beta_1-3/2-n/2}\partial_x^{n+3}v_m \|_{L_{xyt_1}^2}
                \\&
                \leq c(T) \| v_m \|_{L_T^\infty H_{xy}^{n+3}}^2 \| \langle x \rangle^{\frac{n+3-k}{{n+3}}}\partial_x^k v_m \|_{L_T^\infty L_{xy}^2} \| \langle x \rangle^{\beta_1-3/2-n/2}\partial_x^{n+3}v_m \|_{L_{xyt_1}^2}
                \\&
                \leq c(T,\eta_{n+3}(v),\|v\|_{Z_{n+3,\frac{n+3}{2}}}),
            \end{split}
        \end{equation}
        where we have used \eqref{A3} and \eqref{A16}.

\noindent\underline{If $j=n+3$}:\\

In case $k=0$ or $k=n+3$ the estimate is the exact same estimate  \eqref{A21}. In case $k=1$ or $k=n+2$ we have    
\begin{equation}
        \label{A25}
            \begin{split}
                \| (\partial_x v_m)^2 &\partial_x^{n+2}v_m \partial_x^{n+3}v_m \rho_N \|_{L_{xyt_1}^1}\\&
                \leq \| \partial_xv_m \|_{L_{xyt_1}^\infty}^2 
                \| \langle x \rangle^{1/2}\partial_x^{n+2}v_m \|_{L_{xyt_1}^2}
                \| \langle x \rangle ^{\beta_1-3/2-n/2} \partial_x^{n+3}v_m \|_{L_{xyt_1}^2}\\&\leq c(T,\eta_{n+3}(v),\|v\|_{Z_{n+3,\frac{n+3}{2}}}),
            \end{split}
        \end{equation}
        where we have used the inductive hypotheses and \eqref{A16}.\\ 
        
Also, in case $k\in \{ 2,3,\dots, n+1 \}$ we argue as follows
\begin{equation}
            \label{A26}
            \begin{split}
                \| \partial_x v_m &\partial_x^{n+3-k}v_m \partial_x^k v_m \partial_x^{n+3}v_m \rho_N \|_{L_{xyt_1}^1}
                \\ &
                \leq \| \partial_x^{n+3-k}v_m \|_{L_{xyt_1}^\infty} \| \partial_x^k v_m \|_{L_{xyt_1}^\infty} \| \langle x \rangle^{1/2}\partial_x v_m\|_{L_{xyt_1}^2} \| 
                \langle x \rangle^{\beta_1-3/2-n/2 }\partial_x^{n+3}v_m \|_{L_{xyt_1}^2} 
                \\&
                \leq c(T) \| v_m \|_{L_T^\infty H_{xy}^{n+3}}^2 \| \langle x \rangle^{\frac{n+3-1}{2}}\partial_x v_m \|_{L_T^\infty L_{xy}^2} \| \langle x \rangle^{\beta-3/2-n/2}\partial_x^{n+3}v_m \|_{L_{xyt_1}^2}                 \\&
                \leq c(T,\eta_{n+3}(v),\|v\|_{Z_{n+3,\frac{n+3}{2}}}),
            \end{split}
        \end{equation}
        where we have used \eqref{A3} and \eqref{A16}.\\

We now proceed to estimate the $\bar{B}_{jk}$ terms.\\

\noindent\underline{If $j=0$}: Analogous to \eqref{A20} with $\partial_x^{n+3}\partial_y v_m$ instead of $\partial_x^{n+4} v_m$.\\

\noindent\underline{If $j=1$}:\\
\begin{equation}
            \begin{split}
                \| \partial_x^{n+2}&\partial_y v_m \partial_x v_m v_m \partial_x^{n+3} v_m \rho_N \|_{L_{xyt_1}^1}
                \\&
                \leq c(T) \| \partial_x v_m \|_{L_{xyt_1}^\infty} \| \langle x \rangle^{1/2} v_m \|_{L_T^3L_{xy}^\infty} \| \partial_x^{n+2}\partial_y v_m \|_{L_T^\infty L_{xy}^2 } \| \langle x \rangle^{\beta_1-\frac{n+3}{2}}\partial_x^{n+3}v_m \|_{L_{xyt_1}^2}
                \\&
                \leq c(T) \| v_m \|_{L_T^\infty H_{xy}^{n+3}}^2 \| \langle x \rangle^{{\frac{9}{8}}^+} v_m \|_{L_T^3L_{xy}^\infty} \| \langle x \rangle^{\beta_1 -3/2-n/2}\partial_x^{n+3}v_m \|_{L_{xyt_1}^2}
                \\&
                \leq c(T,\eta_{n+3}(v),\|v\|_{Z_{n+3,\frac{n+3}{2}}}),
            \end{split}
        \end{equation}
        where we used Lemma \ref{remark2} and \eqref{A16}.\\

\noindent\underline{If $j\in\{2,3,\dots, n+1 \}$}: Analogous to \eqref{A22} with $\partial_x^{n+3-j}\partial_y v_m$ instead of $\partial_x^{n+4-j}v_m$considering the interpolation 
    \begin{equation*}
        \| \langle x \rangle^{1/2}\partial_x^{n+3-j}\partial_y v_m \|_{L_{xy}^2}\leq 
        c\| J_x^{n+2}\partial_yv_m \|_{L_{xy}^2}^{1-\theta} \| \langle x \rangle^{\frac{n+2}{2(j-1)}} \partial_y v_m \|_{L_{xy}^2}^\theta
    \end{equation*}
    with $\theta= \displaystyle \frac{j-1}{n+2}$, $a=n+2$ and $b=\displaystyle \frac{n+2}{2(j-1)}$ and 
    \begin{align}        
    \label{A28}
        \begin{split}
            \| \langle x \rangle^{\frac{n+2}{2(j-1)}}\partial_y v_m \|_{L_{xy}^2} &\leq
            \| J_{xy}^1 \langle x \rangle^{\frac{n+2}{2(j-1)}}v_m \|_{L_{xy}^2}\\
            &\leq c \| J_{xy}^{n+3}v_m \|_{L_{xy}^2}^\lambda \| \langle x \rangle^{\frac{n+3}{2}} v_m \|_{L_{xy}^2}^{1-\lambda}
        \end{split}
    \end{align}
    with $\lambda=\frac{1}{n+3}$, $a=n+3$, $b=\frac{n+3}{2(j-1)}\leq \frac{n+3}{2}$.\\

\noindent\underline{If $j=n+2$}: Analogous to \eqref{A23} and \eqref{A24} with $\partial_x\partial_y v_m$ instead of $\partial_x^2 v_m$.\\

\noindent\underline{If $j=n+3$}: Analogous to \eqref{A21}, \eqref{A25}, and \eqref{A26} with $\partial_y v_m$ instead of $\partial_x v_m$, where in \eqref{A26} we interpolate as follows:
    \begin{equation}
        \label{A29}
        \| \langle x \rangle^{1/2} \partial_y v_m \|_{L_T^\infty L_{xy}^2}\leq \| J_{xy}^1\langle x \rangle^{1/2} v_m \|_{L_T^\infty L_{xy}^2}\leq c
        \| J_{xy}^{n+3}v_m \|_{L_T^\infty L_{xy}^2}^\theta
        \| \langle x \rangle^{\frac{n+3}{2}}v_m \|_{L_T^\infty L_{xy}^2}^{1-\theta},
    \end{equation}
with $\theta = \frac{n+2}{n+3}$, $a=n+3$ and $b=\frac{n+3}{2(n+2)}\leq \frac{n+3}{2}$.\\

Taking into account \eqref{A20}-\eqref{A29}, integrating \eqref{A19} over $[0,t_1]$ and arguing as done in the proof of Theorem \ref{teo1} it can be seen that 
\begin{equation*}
    \int_0^{t_1}\int (\partial_x^{n+4}v)^2\langle x \rangle^{\beta_1-2-n/2}dxdydt<\infty,
\end{equation*}
which implies 
\begin{equation}
    \label{A30}
    \partial_x^{n+4} v(t)\langle x \rangle^{\beta_1-2-n/2}\in L^2(\R^2) \text{ a.e. }t\in [0,t_1].
\end{equation}
Finally, $P(n+1)$ follows from the choice of $\beta_1$, \eqref{A17}, \eqref{A30} and symmetry of the arguments.
Theorem \ref{T4.2} is proven.
\end{proof}

It is clear that Theorems \ref{teo1}, and \ref{T4.2} imply the results stated in Theorem \ref{1.1}.

\section*{Acknowledgments}
E.B. and J.J.U. were partially supported by Universidad Nacional de Colombia, Sede-Medellín– Facultad de Ciencias – Departamento de Matemáticas – Grupo de investigación en Matemáticas de la Universidad Nacional de Colombia Sede Medellín, carrera 65 No. 59A -110, post 50034, Medellín Colombia. Proyecto: Análisis no lineal aplicado a problemas mixtos en ecuaciones diferenciales parciales, código Hermes 60827. Fondo de Investigación de la Facultad de Ciencias empresa 3062.\\

A.M. is partially supported by CNPq, Conselho Nacional de Desenvolvimento Científico e Tecnológico - Brazil grant 170246/2023-0 and INCTMat, Instituto Nacional de Ciência e Tecnologia de Matemática - Brazil.\\


\end{document}